\documentclass[12pt]{amsart}
\usepackage{amsmath}
\usepackage{amsfonts}
\usepackage{amssymb}
\usepackage{bbding}
\usepackage{stmaryrd}
\usepackage{txfonts}
\usepackage{graphicx}
\usepackage{epsfig}
\usepackage{xypic}
\usepackage{tikz}
\usepackage{longtable}
\usepackage[all]{xy}
\usepackage{pgflibraryarrows}
\usepackage{pgflibrarysnakes}
\usepackage[shortlabels]{enumitem}
\usepackage{ifpdf}
\ifpdf
  \usepackage[colorlinks,final,hyperindex]{hyperref}
\else
  \usepackage[colorlinks,final,hyperindex,hypertex]{hyperref}
\fi
\usepackage{xspace}

\newtheorem{theorem}{Theorem}[section]
\newtheorem{lemma}[theorem]{Lemma}
\newtheorem{coro}[theorem]{Corollary}

\newtheorem{prop}[theorem]{Proposition}
\theoremstyle{definition}
\newtheorem{defn}[theorem]{Definition}
\newtheorem{remark}[theorem]{Remark}
\newtheorem{exmp}[theorem]{Example}

\newcommand{\nc}{\newcommand}
\newcommand{\delete}[1]{}

\nc{\name}[1]{{\bf #1}}
\nc{\tred}[1]{\textcolor{red}{#1}}
\nc{\tblue}[1]{\textcolor{blue}{#1}} \nc{\tgreen}[1]{\textcolor{green}{#1}} \nc{\tpurple}[1]{\textcolor{purple}{#1}} \nc{\btred}[1]{\textcolor{red}{\bf #1}} \nc{\btblue}[1]{\textcolor{blue}{\bf #1}} \nc{\btgreen}[1]{\textcolor{green}{\bf #1}} \nc{\btpurple}[1]{\textcolor{purple}{\bf #1}}

\newcommand{\efootnote}[1]{}

\nc{\mlabel}[1]{\label{#1}}  
\nc{\mcite}[2][]{\cite[#1]{#2}}  
\nc{\mref}[1]{\ref{#1}}  
\nc{\mbibitem}[1]{\bibitem{#1}} 

\delete{
\nc{\mlabel}[1]{\label{#1}  
{\hfill \hspace{1cm}{\bf{{\ }\hfill(#1)}}}}
\nc{\mcite}[1]{\cite{#1}}  
\nc{\mref}[1]{\ref{#1}{{\bf{{\ }(#1)}}}}  
\nc{\mbibitem}[1]{\bibitem[\bf #1]{#1}} 
}

\renewcommand\geq{\geqslant}
\renewcommand\leq{\leqslant}

\renewcommand\bar[1]{\overline{#1}}

\nc{\nz}{\varepsilon}
\nc{\Id}{\mathrm{Id}}
\nc{\map}[2]{{#2}^{#1}}
\nc{\gp}{B}

\nc{\Irr}{\mathrm{Irr}}
\nc{\vx}{\sigma} \nc{\vy}{\tau} \nc{\dvx}{\sigma^{(1)}} \nc{\dvy}{\tau^{(1)}} \nc{\done}{\vep} \nc{\mcitep}[1]{\mcite{#1}} \nc{\wt}{\mathrm{wt}} \nc{\bre}[1]{|#1|} \nc{\mapmonoid}{\frakM} \nc{\disjoint}{\frakM'}
\nc{\ncpoly}[1]{\langle #1\rangle}  
\nc{\mapm}[1]{\lfloor\!|{#1}|\!\rfloor}
\nc{\diff}[1]{{}^\NC\{ #1 \}} \nc{\disj}[1]{\{{#1}\}'} \nc{\mdisj}[1]{\frakM'(#1)} \nc{\brho}{\bar{\rho}} \nc{\om}{\bar{\frakm}} \nc{\frakn}{\mathfrak n} \nc{\ddeg}[1]{^{(#1)}} \nc{\opset}{X} \nc{\genset}{{Z}} \nc{\NC}{\mathrm{{NC}}} \nc{\leaf}{\mathrm{leaf}} \nc{\twig}{\mathrm{twig}} \nc{\fe}{\mathrm{fl}} \nc{\munderline}[1]{#1} \nc{\bo}{o} \nc{\dep}{\mathrm{depth}} \nc{\ofe}{\mathrm{ofl}} \nc{\dfe}{\mathrm{dfe}} \nc{\fex}{\mathrm{fex}} \nc{\dl}{\mathrm{dlex}} \nc{\db}{\mathrm{db}} \nc{\lex}{\mathrm{lex}} \nc{\clex}{\mathrm{clex}} \nc{\dgp}{\mathrm{dgp}} \nc{\dgx}{\mathrm{dgx}} \nc{\br}{\mathrm{br}} \nc{\obd}{\mathrm{odb}} \nc{\ob}{\mathrm{ob}}
\nc{\pie}{\mathrm{PIE}}
\nc{\rbo}{\mathrm{RBO}}
\nc{\supp}{\mathcal{S}}
\nc{\nul}{\mathcal{Z}}

\nc{\bin}[2]{ (_{\stackrel{\scs{#1}}{\scs{#2}}})}  
\nc{\binc}[2]{ \left (\!\! \begin{array}{c} \scs{#1}\\
    \scs{#2} \end{array}\!\! \right )}  
\nc{\bincc}[2]{  \left ( {\scs{#1} \atop
    \vspace{-1cm}\scs{#2}} \right )}  
\nc{\bs}{\bar{S}} \nc{\cosum}{\sqsubset} \nc{\la}{\longrightarrow} \nc{\rar}{\rightarrow} \nc{\dar}{\downarrow} \nc{\dprod}{**} \nc{\dap}[1]{\downarrow \rlap{$\scriptstyle{#1}$}} \nc{\md}[1]{\bar{#1}} \nc{\uap}[1]{\uparrow \rlap{$\scriptstyle{#1}$}} \nc{\defeq}{\stackrel{\rm def}{=}} \nc{\disp}[1]{\displaystyle{#1}} \nc{\dotcup}{\ \displaystyle{\bigcup^\bullet}\ } \nc{\gzeta}{\bar{\zeta}} \nc{\hcm}{\ \hat{,}\ } \nc{\hts}{\hat{\otimes}} \nc{\barot}{{\otimes}} \nc{\free}[1]{\bar{#1}} \nc{\uni}[1]{\tilde{#1}} \nc{\hcirc}{\hat{\circ}} \nc{\leng}{\ell} \nc{\lleft}{[} \nc{\lright}{]} \nc{\lc}{\lfloor} \nc{\rc}{\rfloor}
\nc{\lb}{[} 
\nc{\rb}{]} 
\nc{\curlyl}{\left \{ \begin{array}{c} {} \\ {} \end{array}
    \right.  \!\!\!\!\!\!\!}
\nc{\curlyr}{ \!\!\!\!\!\!\!
    \left. \begin{array}{c} {} \\ {} \end{array}
    \right \} }
\nc{\longmid}{\left | \begin{array}{c} {} \\ {} \end{array}
    \right. \!\!\!\!\!\!\!}
\nc{\onetree}{\bullet} \nc{\ora}[1]{\stackrel{#1}{\rar}}
\nc{\ola}[1]{\stackrel{#1}{\la}}
\nc{\ot}{\otimes} \nc{\mot}{{{\boxtimes\,}}} \nc{\otm}{\overline{\boxtimes}} \nc{\sprod}{\bullet} \nc{\scs}[1]{\scriptstyle{#1}} \nc{\mrm}[1]{{\rm #1}} \nc{\msum}{\sum\limits}
\nc{\margin}[1]{\marginpar{\rm #1}}   
\nc{\dirlim}{\displaystyle{\lim_{\longrightarrow}}\,} \nc{\invlim}{\displaystyle{\lim_{\longleftarrow}}\,} \nc{\mvp}{\vspace{0.3cm}} \nc{\tk}{^{(k)}} \nc{\tp}{^\prime} \nc{\ttp}{^{\prime\prime}} \nc{\svp}{\vspace{2cm}} \nc{\vp}{\vspace{8cm}} \nc{\proofbegin}{\noindent{\bf Proof: }}
\nc{\proofend}{$\blacksquare$ \vspace{0.3cm}}
\nc{\modg}[1]{\!<\!\!{#1}\!\!>}
\nc{\intg}[1]{F_C(#1)} \nc{\lmodg}{\!<\!\!} \nc{\rmodg}{\!\!>\!} \nc{\cpi}{\widehat{\Pi}}
\nc{\sha}{{\,\mbox{\cyr X}\,}}  
\nc{\shap}{{\mbox{\cyrs X}}} 
\nc{\shpr}{\diamond}    
\nc{\shp}{\ast} \nc{\shplus}{\shpr^+}
\nc{\shprc}{\shpr_c}    
\nc{\msh}{\ast} \nc{\zprod}{m_0} \nc{\oprod}{m_1} \nc{\vep}{\varepsilon} \nc{\labs}{\mid\!} \nc{\rabs}{\!\mid}
\nc{\astarrow}{\overset{\raisebox{-3pt}{$\ast$}}{\rightarrow}}

\nc{\Sym}{\mrm{Sym}}
\nc{\Nsym}{\mrm{NSym}}
\nc{\QSym}{\mrm{QSym}}
\nc{\BQSym}{\mrm{BQSym}}
\nc{\ZQsym}{\mrm{ZQSym}}
\nc{\SSym}{\mrm{\mathfrak{S}Sym}}
\nc{\HSym}{\mrm{\mathfrak{H}Sym}}
\nc{\HHSym}{\mrm{HSym}}
\nc{\syms}{symmetric functions\xspace}
\nc{\qsyms}{quasi-symmetric functions\xspace}
\nc{\nsymg}{\mathrm{NSym}_\gp}
\nc{\allow}{{\rm {allow}}}
\nc{\Ap}{{\rm {Ap}}}
\nc{\ass}{{\rm{ass}}}
\nc{\cmpt}{{\rm{cmpt}}}
\nc{\cmpts}{{\rm{cmpts}}}
\nc{\Comp}{{\rm {Comp}}}
\nc{\des}{{\rm {des}}}
\nc{\Des}{{\rm {Des}}}
\nc{\dom}{{\rm{dom}}}
\nc{\GDes}{{\rm {GDes}}}
\nc{\Inv}{{\rm{Inv}}}
\nc{\inv}{{\rm{inv}}}
\nc{\Invol}{{\rm{Invol}}}
\nc{\maj}{{\rm{maj}}}
\nc{\Nega}{{\rm{Neg}}}
\nc{\nega}{{\rm{neg}}}
\nc{\Nsp}{{\rm{Nsp}}}
\nc{\nsp}{{\rm{nsp}}}
\nc{\parr}{{\rm {Par}}}
\nc{\PC}{{\mathcal{PC}}}
\nc{\Pos}{{\rm{Pos}}}
\nc{\ran}{{\rm{ran}}}
\nc{\rank}{{\rm{rank}}}
\nc{\Set}{{\rm {Set}}}
\nc{\set}{{\rm {set}}}
\nc{\Sh}{{\rm {Sh}}}
\nc{\sign}{{\rm {sign}}}
\nc{\Span}{{\rm {Span}}}
\nc{\st}{{\rm{st}}}
\nc{\stb}{{\rm{stb}}}
\nc{\sts}{{\rm{sts}}}
\nc{\wcomp}{{\rm {wcomp}}}
\nc{\wcts}{|\hspace{-1.2mm}\models}
\nc{\Sol}{\mrm {Sol}}

\nc{\dth}{d} \nc{\mmbox}[1]{\mbox{\ #1\ }} \nc{\fp}{\mrm{FP}} \nc{\rchar}{\mrm{char}} \nc{\Fil}{\mrm{Fil}} \nc{\Mor}{Mor\xspace} \nc{\gmzvs}{gMZV\xspace} \nc{\gmzv}{gMZV\xspace} \nc{\mzv}{MZV\xspace} \nc{\mzvs}{MZVs\xspace} \nc{\Hom}{\mrm{Hom}} \nc{\id}{\mrm{id}} \nc{\im}{\mrm{im}} \nc{\incl}{\mrm{incl}}  \nc{\mchar}{\rm char}

\nc{\Alg}{\mathbf{Alg}} \nc{\Bax}{\mathbf{Bax}} \nc{\bff}{\mathbf f} \nc{\bfk}{{\bf k}} \nc{\bfone}{{\bf 1}} \nc{\bfx}{\mathbf x} \nc{\bfy}{\mathbf y}
\nc{\base}[1]{\bfone^{\otimes ({#1}+1)}} 
\nc{\Cat}{\mathbf{Cat}} \delete{}
\nc{\detail}{\marginpar{\bf More detail}
    \noindent{\bf Need more detail!}
    \svp}
\nc{\Int}{\mathbf{Int}} \nc{\Mon}{\mathbf{Mon}}
\nc{\rbtm}{{shuffle }} \nc{\rbto}{{Rota-Baxter }} \nc{\remarks}{\noindent{\bf Remarks: }} \nc{\Rings}{\mathbf{Rings}} \nc{\Sets}{\mathbf{Sets}}
\nc{\balpha}{\mathbf{\alpha}}

\nc{\BA}{{\mathbb A}} \nc{\CC}{{\mathbb C}} \nc{\DD}{{\mathbb D}} \nc{\EE}{{\mathbb E}} \nc{\FF}{{\mathbb F}} \nc{\GG}{{\mathbb G}} \nc{\HH}{{\mathbb H}} \nc{\LL}{{\mathbb L}} \nc{\NN}{{\mathbb N}} \nc{\KK}{{\mathbb K}} \nc{\PP}{{\mathbb P}} \nc{\QQ}{{\mathbb Q}} \nc{\RR}{{\mathbb R}} \nc{\TT}{{\mathbb T}} \nc{\VV}{{\mathbb V}} \nc{\ZZ}{{\mathbb Z}}

\nc{\cala}{{\mathcal A}} \nc{\calc}{{\mathcal C}} \nc{\cald}{{\mathcal D}} \nc{\cale}{{\mathcal E}} \nc{\calf}{{\mathcal F}} \nc{\calg}{{\mathcal G}} \nc{\calh}{{\mathcal H}} \nc{\cali}{{\mathcal I}} \nc{\call}{{\mathcal L}} \nc{\calm}{{\mathcal M}} \nc{\caln}{{\mathcal N}} \nc{\calo}{{\mathcal O}} \nc{\calp}{{\mathcal P}} \nc{\calr}{{\mathcal R}} \nc{\cals}{{\mathcal S}} \nc{\calt}{{\mathcal T}} \nc{\calw}{{\mathcal W}} \nc{\calk}{{\mathcal K}} \nc{\calx}{{\mathcal X}}
\nc{\calz}{{\mathcal Z}}
 \nc{\CA}{\mathcal{A}}

\nc{\fraka}{{\mathfrak a}} \nc{\frakA}{{\mathfrak A}} \nc{\frakb}{{\mathfrak b}} \nc{\frakB}{{\mathfrak B}}
\nc{\frakc}{{\mathfrak c}}  \nc{\frakD}{{\mathfrak D}}
\nc{\frakH}{{\mathfrak H}}
\nc{\frakh}{{\mathfrak h}} \nc{\frakM}{{\mathfrak M}}
\nc{\frakO}{{\mathfrak O}}
\nc{\frakE}{{\mathfrak E}}
\nc{\bfrakM}{\overline{\frakM}} \nc{\frakm}{{\mathfrak m}} \nc{\frakP}{{\mathfrak P}} \nc{\frakN}{{\mathfrak N}} \nc{\frakp}{{\mathfrak p}} \nc{\frakS}{{\mathfrak S}}
\nc{\frakk}{{\mathfrak k}}
\nc{\frakx}{{\mathfrak x}}
\nc{\frakl}{{\mathfrak l}} \nc{\ox}{\bar{\frakx}} \nc{\frakX}{{\mathfrak X}} \nc{\fraky}{{\mathfrak y}} \nc\dop{\delta}
\nc{\Reduce}{{\rm Red}}

\font\cyr=wncyr10 \font\cyrs=wncyr7
\nc{\redt}[1]{\textcolor{red}{#1}}
\nc{\yu}[1]{\textcolor{red}{\tt Yu:#1}}

\begin{document}
\title[Weak order and the monomial basis of $\HSym$]
{The weak order on the hyperoctahedral group and the monomial basis for the Hopf algebra of signed permutations}

\author{Houyi Yu}
\address{School of Mathematics and Statistics, Southwest University, Chongqing 400715, China}
\email{yuhouyi@swu.edu.cn}

\hyphenpenalty=8000
\date{}

\begin{abstract}
We give a combinatorial description for the weak order on the hyperoctahedral group.
This characterization is then used to analyze the order-theoretic properties of the shifted products of hyperoctahedral groups.
It is shown that each shifted product is a disjoint union of some intervals, which can be convex embedded into a hyperoctahedral group.
As an application, we investigate the monomial basis for the Hopf algebra $\HSym$ of signed permutations, related to the fundamental basis via M\"obius inversion on the weak order on hyperoctahedral groups.
It turns out that the image of a monomial basis element under the descent map from $\HSym$ to the algebra of type $B$ quasi-symmetric functions is either zero or a monomial quasi-symmetric function of type $B$.
\end{abstract}

\keywords{ weak order, hyperoctahedral group, Hopf algebra of signed permutations, monomial basis}

\maketitle

\tableofcontents

\hyphenpenalty=8000 \setcounter{section}{0}


\allowdisplaybreaks

\section{Introduction}\label{sec:int}

The (left) weak order is a powerful tool in the combinatorial study of a Coxeter group.
It can be defined as the suffix order of reduced expressions of the group elements,
or more combinatorially as the inclusion order of the right associated reflection sets of the group elements.
The general structure of the weak order on an arbitrary Coxeter group was first systematically investigated by Bjorner \cite{Bj84},
and extensively studied in various recent works, including \cite{BB03,CPP19,Dhp18,Dy19,GP20,Re19}.

The Coxeter groups of types $A$ and $B$ have combinatorial interpretations as the symmetric group $\mathfrak{S}_n$ of permutations and
the hyperoctahedral group $\mathfrak{B}_n$ of signed permutations, respectively.
The weak order on $\mathfrak{S}_n$ was first studied by Yanagimoto and Okamoto \cite{Ya69}, who showed that
it coincides with the inclusion order of the inversion sets of permutations, where an \emph{inversion} of a permutation $w\in\mathfrak{S}_n$
is a pair $(i,j)$ such that $1\leq i<j\leq n$ and $w(i)>w(j)$.
Somewhat surprisingly, there is not an analogous description for the weak order on $\mathfrak{B}_n$.
One of the aims of this paper is to provide a combinatorial characterization for the weak order on $\mathfrak{B}_n$, and apply this description to investigate the order structures of the shifted products of hyperoctahedral groups.

The weak order is not only of interest for their own in combinatorics but also closely related to the structures of various algebras arising from Coxeter groups.
Given positive integers $p_1,p_2,\ldots,p_k$, where $k\geq 2$,
for any $(p_1,p_2,\ldots,p_k)$-shuffle $\xi$, Aguiar and Sottile \cite[Proposition 2.10]{AS05} showed that the map
\begin{align}\label{eq:spqxi}
    \mathfrak{S}_{p_1}\times \mathfrak{S}_{p_2}\times \cdots\times  \mathfrak{S}_{p_k}\hookrightarrow\mathfrak{S}_{p_1+p_2+\cdots+p_k},\quad u\mapsto u \xi^{-1}
\end{align}
is a convex embedding under the weak order on $\mathfrak{S}_{p_1+p_2+\cdots+p_k}$, so that it preserves meets and joins.
When $k=2$, this leads to a connected graded self-dual Hopf algebra structure on the space $\SSym:=\bigoplus_{n\geq0}\bfk \mathfrak{S}_n$, which has
a basis $\{F_u\,|\,u\in \mathfrak{S}_n,n\geq0\}$ called the \emph{fundamental basis}.
It is the well-known \emph{Malvenuto-Reutenauer Hopf algebra} of permutations \cite{MR95},
whose algebraic structure is well understood in terms of the weak order on the symmetric group \cite{AS05,LR02}.

In order to analyze the detailed Hopf structure of $\SSym$, Aguiar and Sottile \cite{AS05} produced a new basis, called the \emph{monomial basis}, which is related to the fundamental basis by M\"obius inversion on the weak order on the symmetric group, and provided enumerative-combinatorial descriptions of the Hopf algebra operations in terms of this basis. These results are all based on the convexities of the embeddings defined by Eq.\,\eqref{eq:spqxi}.

Recently, the techniques of Aguiar and Sottile \cite{AS05} are abstracted by Bergeron, D'le\'on, Li, Pang and Vargas \cite{BDLPV21} to form a set of axioms so that one can define a monomial basis on any combinatorial Hopf algebra. These axioms guarantee that the monomial basis enjoys many remarkable algebraic properties, including a positive multiplication formula and a cancellation-free antipode formula.

There is an analogous construction for hyperoctahedral groups.
Given a $(p_1,p_2,\ldots,p_k)$-shuffle $\xi$, we have a map corresponding to $\xi$ defined by
\begin{align}\label{eq:spqxib}
    \mathfrak{B}_{p_1}\times \mathfrak{B}_{p_2}\times \cdots\times  \mathfrak{B}_{p_k}\hookrightarrow\mathfrak{B}_{p_1+p_2+\cdots+p_k},\quad u\mapsto u \xi^{-1},
\end{align}
which induces a connected graded Hopf algebra structure on the space $\HSym:= \bigoplus_{n\geq0}\bfk \mathfrak{B}_n$,
whose $n$th graded component has the \emph{fundamental basis} $\{\mathbf{F}_u\,|\,u\in \mathfrak{B}_n\}$.
This Hopf algebra structure with shifted shuffle product was first considered by Aguiar, Bergeron and Nyman \cite{ABN04},
and later in more general setting by Novelli and Thibon \cite{NT10}.
Quite recently a generalization of this structure was introduced in \cite{GTY20}, where the product is replaced by the shifted quasi-shuffle product with a weight, leading to a new
Hopf algebra beyond the category of combinatorial Hopf algebras in the sense of \cite{ABS06}, since
the product is in general no longer graded with respect to the size of signed permutations.

Mainly motivated by the work in \cite{AS05,BDLPV21},
we investigate the monomial basis $\{\mathbf{M}_{u}\,|\,u\in \mathfrak{B}_n,n\geq0\}$ for the Hopf algebra $\HSym$ with respect to the weak order on $\mathfrak{B}_n$,
where,  for $u\in\mathfrak{B}_n$,
\begin{align*}
    \mathbf{M}_{u}=\sum_{u\leq v}\mu_{\mathfrak{B}_n}(u,v)\mathbf{F}_{v},\quad \text{or equivalently,}\quad \mathbf{F}_{u}=\sum_{u\leq v}\mathbf{M}_{v}.
\end{align*}
Here $\mu_{\mathfrak{B}_n}$ is the M\"obius function of the weak order on $\mathfrak{B}_n$.
To this end, we analyze in detail the order structure of $\mathfrak{B}_{p_1}\times \mathfrak{B}_{p_2}\times \cdots\times  \mathfrak{B}_{p_k}$.

It turns out that the  Hopf algebra $\HSym$ under the weak order on hyperoctahedral groups does not satisfy the axioms established in \cite{BDLPV21}.
However, the monomial basis $\mathbf{M}_{u}$ is well-behaved in the sense that it is compatible with the algebra $\BQSym$ of quasi-symmetric functions of type $B$.
Chow \cite{Ch01} showed that there is an algebra homomorphism from $\HSym$ onto $\BQSym$ that sends one fundamental basis to the other, induced by taking descent sets of signed permutations.
We obtain an explicit formula for this map in terms of the monomial basis.
Our main tool is the combinatorial description of the weak order on $\mathfrak{B}_n$.

The structure of the paper is as follows. After summarizing the notation and basic facts on Coxeter groups in Section \ref{sec:woocg}, we give in
Section \ref{sec:lwoohg} a combinatorial characterization of the weak order on the hyperoctahedral group.
This description is an effective criteria which enables us to compare directly two signed permutations of the same size in the weak order by verifying their inversion sets, negative index sets and negative sum pair sets, respectively.
Section \ref{sec:opmonC} is devoted to studying the order structure of the shifted products of hyperoctahedral groups.
Aguiar and Sottile \cite[Proposition 2.10]{AS05}
showed that any parabolic subgroup of a symmetric group can be convex embedded into the symmetric group.
We first generalize this result from type $A$ to all Coxeter groups, and then study the order structure of the shifted products of hyperoctahedral groups.
In particular, we show that each component of a shifted product of hyperoctahedral groups is an interval isomorphic to the unique component which is a parabolic subgroup. Moreover, all components can be convex embedded into a hyperoctahedral group under the restriction of the map defined by Eq.\,\eqref{eq:spqxib}.
These results are applied in Section \ref{sec:mbhsp} to study the monomial basis of the Hopf algebra $\HSym$ of signed permutations.
We give sufficient conditions for the structure constants of the coproduct and product in terms of the monomial basis to be nonnegative.
It is shown that the image of a monomial basis element under the descent map from $\HSym$ to $\BQSym$ is either zero or a monomial quasi-symmetric function of type $B$.

\section{Preliminaries}\label{sec:woocg}

Unless otherwise specified,
for any nonnegative integers $m,n$ with $m\leq n$, let $[m,n]=\{m,m+1,\ldots,n\}$ and $[n]=[1,n]$ if $n\geq1$.
For a set $I$ of integers, the notation $I = \{i_1 < i_2 < \cdots < i_k \}$ indicates that $I = \{i_1, i_2, \cdots, i_k \}$ and
$i_1 < i_2 < \cdots < i_k$. The cardinality of a finite set $A$ will be denoted by $\#A$.

\subsection{Coxeter groups}\label{subsec:lwooq}
We begin by recalling some necessary background results on posets and Coxeter groups.
For a detailed treatment we refer to \cite{BB05,Bb81,Hum92,Sta12}.

Let $(P,\leq_P)$ be a poset. When no possible
confusion may arise, we will simply denote $\leq_P$ by $\leq$ and write $P$ for $(P,\leq)$.
If $x,y\in P$, then $y$ \emph{covers} $x$ or $x$ is \emph{covered} by $y$, denoted $x\prec y$ or $y\succ x$,
if $x< y$ and there is no element $z\in P$ such that $x< z< y$.
For $x\leq y$ in $P$, the \emph{closed interval} $[x,y]$ is the subset $\{z\in P\,|\,x\leq z\leq y\}$, endowed with the partial order induced from $P$.

A \emph{Coxeter system} is a pair $(W,S)$, where $W$ is a group and $S$ is a set of generators of $W$ subject to the relations
$$
(ss')^{m(s,s')}=e \ \text{for all} \ s,s'\in S,
$$
where $m(s,s')$ denotes the order of $ss'$ and $m(s,s')=1$ if and only if $s=s'$. The group $W$ is called  a \emph{Coxeter group} and $S$ is the set of \emph{Coxeter generators}.


Let $(W,S)$ be a Coxeter system. For any element $w$ of $W$,
the \emph{length} of $w$, denoted $\ell(w)$, is the least $k$ such that
$w=s_1s_2\cdots s_k$ with all $s_i\in S$. Such a decomposition
is called a \emph{reduced word} (or \emph{reduced expression}) for $w$.
Let $T:=\{wsw^{-1}\,|\,s\in S, w\in W\}$ be the set of \emph{reflections} of $W$.
Then $\ell(w)=\#T_R(w)$ for any $w\in W$, where
$
T_R(w):=\{t\in T\,|\, \ell(wt)<\ell(w)\}.
$

For $u$ and $v$ in a Coxeter group $W$, we say that $u$ precedes $v$ in \emph{left weak order}, written $u\leq v$, if $\ell(v)=\ell(u)+\ell(v u^{-1})$.
The right weak order is defined similarly, with $u^{-1}v$ instead of $vu^{-1}$.
We will only use the left weak order and refer to it as the weak order for short. We abuse notation and identify the Coxeter group $W$ with the poset consisting of $W$ equipped
with the weak order.
Any finite Coxeter group $W$ is a graded lattice with the least element $e$. In particular, the rank function is given by length.
For later reference, we state the following well-known useful characterizations of the weak order, which can be found in \cite{BB05}.

\begin{lemma}\label{lem:uleqlvequiv}
Let $(W,S)$ be a Coxeter system, and $u,v\in W$. Then the following conditions are equivalent:
\begin{enumerate}
\item\label{lemitemleqLa} $u\leq v$;
\item\label{lemitemleqLb} $T_R(u)\subseteq T_R(v)$;
\item\label{lemitemleqLd} {\rm (Suffix Property)} there exist reduced words $u=s_1s_2\cdots s_k$ and $w=s_1's_2'\cdots s_q's_1s_2\cdots s_k$.
\end{enumerate}
\end{lemma}

Given a subset $J$ of $S$,
the subgroup $W_J$ generated by $J$ is called a \emph{parabolic subgroup} of $W$. Hence, $W_{\emptyset}=\{e\}$ and $W_S=W$.
It is useful to mention that $(W_J,J)$ is also a Coxeter group and the length function of $W_J$ with respect to $J$ coincides with length function of $W$ with respect to $S$.
The same is true for the weak order.
Fix a subset $J$ of $S$, then
\begin{align}\label{eq:w^Jset}
W^J:=\{w\in W\,|\,\ell(w)<\ell(ws)\ \text{for all}\ s\in J\}
\end{align}
is a complete set of the left coset representatives of $W_J$ in $W$, consisting of the unique representatives
of minimal length. 
When $W$ is a finite Coxeter group, $W_J$ and $W^J$ are both intervals \cite{BW88}.
One basic fact about $W^J$ is that every $w\in W$ has a uniquely factorization $w=w^J\cdot w_J$ such that $w^J\in W^J$ and $w_J\in W_J$. Moreover, this factorization satisfies  $\ell(w)=\ell(w^J)+\ell(w_J)$.
The map $w\mapsto w_J$ is a lattice homomorphism from $W$ to $W_J$, which is called a \emph{parabolic homomorphism},
while the map $w\mapsto w^J$ is order-preserving from $W$ to $W^J$ but is not in general a lattice homomorphism \cite{Je05,Re04}.


\subsection{Coxeter groups of types $A$ and $B$}
We refer to \cite[Section 8.1]{BB05} for details on signed permutations, and review some related
notation here.

For a positive integer $n$, we write $[\pm n]$ for the set $\{0,\pm 1,\pm 2,\ldots,\pm n\}$, and take the natural order of integers on $[\pm n]$.
For ease of notation, we use the bar to denote a negative sign, so $\bar{i}=-i$ for $i\in[n]$.
A \emph{signed permutation of size $n$} is a permutation $w$ on $[\pm n]$ satisfying $w(\bar{i})=\bar{w(i)}$.
Notice that $w(0)=0$ and the element $w$ is completely determined by $w(1),w(2), \cdots, w(n)$. So, in one-line notation, we write
$w=(w_1,w_2,\ldots, w_n)$, or simply denote $w=w_1w_2\cdots w_n$,
where $w_i=w(i)$ for $i\in [n]$.
For example, $w=2\bar5 1 \bar3\ \bar4$ is a signed permutation on $[\pm5]$ with  $w(2)=\bar5$,  $w(3)=1$ and $w(\bar{4})=3$.
The set of all signed permutations on $[\pm n]$ naturally form a group under composition, called the \emph{$n$-th hyperoctahedral group} and denoted by $\mathfrak{B}_n$,
which is the Coxeter group of type $B_n$.
The element $1_n$ denotes the identity of $\mathfrak{B}_n$.
Let $s_0$ be the permutation swapping $1$ and $-1$, and for $i\in[n-1]$,
let $s_i$ be the product of transpositions $(i, i+1)(\bar{i},\bar{i+1})$. Write $S_n^B:=\{s_0,s_1,\ldots,s_{n-1}\}$. Then $(\mathfrak{B}_n,S_n^B)$ is a Coxeter system of type $B_n$.

We identify the $n$-th symmetric group $\mathfrak{S}_n$ with the subgroup of $\mathfrak{B}_n$ consisting of all signed permutations $w$ such that $w([n]) =[n]$.
Then $\mathfrak{S}_n$ is the Coxeter group of type $A_{n-1}$ consisting of all permutations of $[n]$, and the set $S_n^A$ of Coxeter generators consists of the elementary transpositions $s_i=(i,i+1)$ for $i\in[n-1]$.
For the rest of this paper, if there is no danger of confusion, we write simply $S^A$ and $S^B$ instead of $S_n^A$ and $S_n^B$, respectively.

Let $w$ be a signed permutation of size $n$.
A positive integer $i\in[n]$ is a \emph{negative index} of $w$ if $w_i<0$.
Given a  pair $(i,j)\in [n]\times [n]$, the pair $(i,j)$ is called an \emph{inversion} of $w$ if $i<j$ and $w_i>w_j$, while $(i,j)$ is called a \emph{negative sum pair} of $w$ if $i<j$ and $w_i+w_j<0$.
Let
$$\nega(w):=\#\Nega(w),\quad
\inv(w):=\#\Inv(w)\quad \text{and} \quad \nsp(w):=\#\Nsp(w),$$ where
\begin{align*}
\Nega(w):=\{i\in[n]\,|\,w_i<0\},\quad \Inv(w):=\{(i,j)\in [n]\times [n]\,|\,i<j,w_i>w_j\}
\end{align*}
and
\begin{align*}
\Nsp(w):=\{(i,j)\in [n]\times [n]\,|\,i<j,w_i+w_j<0\}.
\end{align*}
Then the length of $w$ is given by
\begin{align}\label{eq:ell(w)b}
\ell(w)=\inv(w)+\nega(w)+\nsp(w)=\inv(w)-\sum_{i\in[n], w_i<0}w_i.
\end{align}
In particular, if  $w$ is a permutation, then the length of $w$ is equal to the number of inversions of $w$, that is, $\ell(w)=\inv(w)$.

Let $a=a_1a_2\cdots a_n$ be a word of $n$ nonzero integers. The \emph{standard permutation} $\st(a)$ of $a$ is the unique permutation $w\in\mathfrak{S}_n$ defined by
$$
w_i<w_j\quad \Leftrightarrow\quad a_i\leq a_j
$$
for all $i,j$ with $1\leq i<j\leq n$, while
the \emph{standard signed permutation} $\sts(a)$ of $a$ is the unique signed permutation $w\in \mathfrak{B}_n$ such that $\Nega(w)=\Nega(a)$ and
\begin{align*}
|w_i|<|w_j|\ \Leftrightarrow\ |a_i|\leq |a_j|
\end{align*}
for all $i,j$ with $1\leq i<j\leq n$. Here $\Nega(a):=\{i\in[n]\,|\, a_i<0\}$, and $|m|$ is the absolute value of $m$ for an integer $m$.
One can obtain $\sts(a)$ by first taking absolute values of the numbers, then extracting the standard permutation, and finally putting back the signs of the original numbers.
Clearly, $\st$ coincide with $\sts$ when restrict to words on positive integers. For instance, $\st(352)=\sts(352)=231$, $\st(6\bar{3}2\bar{7}35)=623145$ and $\sts(6\bar{3}2\bar{7}35)=5\bar{2}1\bar{6}34$.

Let $J=S^A\backslash\{s_{p_1},s_{p_2},\ldots,s_{p_k}\}$ where $1\leq p_1<p_2<\cdots< p_k\leq n-1$.
Then the parabolic subgroup of $\mathfrak{S}_n$ generated by $J$ is the subgroup
$$
\mathfrak{S}_J:=\mathfrak{S}_{p_1}\times \mathfrak{S}_{p_2-p_1}\times \cdots\times \mathfrak{S}_{n-p_k},
$$
where $\mathfrak{S}_{p_i-p_{i-1}}$ permutes $[p_{i-1}+1,p_i]$ for $i\in[k+1]$,  with the notation $p_{0}=0$ and $p_{k+1}=n$.
For $u^{(i)}\in \mathfrak{S}_{p_i-p_{i-1}}$ where $i\in [k+1]$, we use $u^{(1)}\times u^{(2)}\times\cdots\times u^{(k+1)}$ to denote the permutation in $\mathfrak{S}_n$ corresponding to
$(u^{(1)}, u^{(2)},\ldots, u^{(k+1)})\in\mathfrak{S}_J$. Then for any $w\in \mathfrak{S}_J$, in terms of the standard permutation of a word, we have
\begin{align*}
    w_J=\st(w_1\cdots w_{p_1})\times \st(w_{p_1+1}\cdots w_{p_2})\times \cdots \times \st(w_{p_k+1}\cdots w_{n}).
\end{align*}
By Eq. \eqref{eq:w^Jset},
the set $\mathfrak{S}^{J}$  of minimal representatives of left cosets of $\mathfrak{S}_{J}$ in $\mathfrak{S}_{n}$ is 
\begin{align*}
\Sh(p_1,p_2-p_1,\ldots,n-p_k)
:=\{w\in \mathfrak{S}_{n}\,|\,w_1<\cdots <w_{p_1},w_{p_1+1}<\cdots <w_{p_2},\ldots,w_{p_k+1}<\cdots <w_{n}\},
\end{align*}
whose elements are usually called \emph{$(p_1,p_2-p_1,\ldots,n-p_k)$-shuffles}.

Now let $J=S^B\backslash\{s_{p_1},s_{p_2},\ldots,s_{p_k}\}$  where $0\leq p_1<p_2<\cdots< p_k\leq n-1$. Then the parabolic subgroup $\mathfrak{B}_J$ has the form
$$
\mathfrak{B}_J:=\mathfrak{B}_{p_1}\times \mathfrak{S}_{p_2-p_1}\cdots\times \mathfrak{S}_{n-p_k}.
$$
Here $p_1=0$ corresponds to the case $s_0\not\in J$ and we write $\mathfrak{B}_{0}=\{1_0\}$, which is the trivial group consisting of only one element.
Hence,
\begin{align*}
\mathfrak{B}_J=\begin{cases}
\mathfrak{S}_{p_2}\times\cdots\times \mathfrak{S}_{n-p_k}, & p_1=0,\\
\mathfrak{B}_{p_1}\times \mathfrak{S}_{p_2-p_1}\cdots\times \mathfrak{S}_{n-p_k}, & p_1>0.
\end{cases}
\end{align*}
 Then for any $w\in \mathfrak{B}_n$, we have
\begin{align}\label{eq:stbtoa}
    w_J=\sts(w_1\cdots w_{p_1})\times \st(w_{p_1+1}\cdots w_{p_2})\times \cdots \times \st(w_{p_k+1}\cdots w_{n}).
\end{align}
The set of minimal left coset representatives of $\mathfrak{B}_J$ is
\begin{align*}
\mathfrak{B}^J=\{w\in \mathfrak{B}_n\,|\, 0<w_1<\cdots<w_{p_1}, w_{p_1+1}<\cdots<w_{p_2},\ldots,w_{p_k+1}<\cdots<w_{n}\}.
\end{align*}

\section{The weak order on the hyperoctahedral group}\label{sec:lwoohg}

Jedli\v{c}ka \cite{Je05} provided a combinatorial construction of the weak order on a Coxeter group by using the semidirect product of semilattices.
In this section, we give a combinatorial description for the weak order on the hyperoctahedral group $\mathfrak{B}_n$
in terms of inversion sets, negative index sets and negative sum pair sets of signed permutations.
This enables us to compare directly two signed permutations when they are expressed in the one-line notation.

A direct translation of the definition of the weak order implies that $u\leq  v$ in $\mathfrak{S}_{n}$ if and only if $\Inv(u)\subseteq \Inv(v)$, while
$u\leq v$ in $\mathfrak{B}_{n}$ if and only if $v=s_iu$ and $u^{-1}(i)<u^{-1}(i+1)$ for some $s_i\in S^B$.
The weak order on $\mathfrak{B}_3$ is illustrated in Figure $1$.

Let $n$ be a positive integer and let $u\leq v$ in $\mathfrak{B}_n$. Note that $\mathfrak{S}_n$ is the parabolic subgroup of $\mathfrak{B}_n$
generated by $S^A$, so it follows from Eq.\,\eqref{eq:stbtoa} that $\st$ is indeed a lattice homomorphism from $\mathfrak{B}_n$
to $\mathfrak{S}_n$ under weak orders, and hence $\st(u)\leq \st(v)$ in $\mathfrak{S}_n$.
Observe that the inversion set of a signed permutation coincides with that of its standard permutation. Thus $\Inv(u)\subseteq\Inv(v)$.
More generally, we have the following result, which can be regarded as a generalization of the weak order on $\mathfrak{S}_n$ and a refinement of
the length function given by Eq.\,\eqref{eq:ell(w)b}.

\begin{theorem}\label{thm:leqliffinvnegnsp}
Let $u,v\in \mathfrak{B}_n$. Then $u\leq v$ if and only if $\Inv(u)\subseteq \Inv(v)$, $\Nega(u)\subseteq\Nega(v)$ and $\Nsp(u)\subseteq\Nsp(v)$.
\end{theorem}

To prove theorem \ref{thm:leqliffinvnegnsp}, we need several lemmas.

\begin{lemma}\label{lem:leqinn}
Let $u,v\in \mathfrak{B}_n$ with $u\prec v$, and let $v=s_iu$ for some $i\in[0,n-1]$.
\begin{enumerate}
  \item If $i=0$, then $\Inv(v)=\Inv(u)$, $\Nsp(v)=\Nsp(u)$ and $\Nega(v)=\Nega(u)\cup\{u^{-1}(1)\}.$
  \item If $i>0$, then $\Nega(v)=\Nega(u)$ and
\begin{align*}
\begin{cases}
\Inv(v)=\Inv(u)\cup\{(u^{-1}(i),\ u^{-1}(i+1))\}, \Nsp(v)=\Nsp(u),& \text{if}\ 0<u^{-1}(i)<u^{-1}(i+1),\\
\Inv(v)=\Inv(u),\ \Nsp(v)=\Nsp(u)\cup\{(u^{-1}(\bar{i}),u^{-1}(i+1))\},& \text{if}\ 0<u^{-1}(\bar{i})<u^{-1}(i+1),\\
\Inv(v)=\Inv(u),\ \Nsp(v)=\Nsp(u)\cup\{(u^{-1}(i+1),u^{-1}(\bar{i}))\},& \text{if}\ 0<u^{-1}(i+1)<u^{-1}(\bar{i}),\\
\Inv(v)=\Inv(u)\cup\{(u^{-1}(\bar{i+1}),u^{-1}(\bar{i}))\},\ \Nsp(v)=\Nsp(u),& \text{if}\ 0<u^{-1}(\bar{i+1})<u^{-1}(\bar{i}).
\end{cases}
\end{align*}
\end{enumerate}
\end{lemma}

\begin{proof}
Assume that $u\prec v$ in $\mathfrak{B}_n$, and let $v=s_iu$ for some $i\in[0,n-1]$.
If $i=0$, then $u\prec s_0u$ so that $u^{-1}(1)>u^{-1}(0)=0$, and hence $v$ is the signed permutation obtained from $u$ by replacing $1$ by $\bar{1}$.
Thus
$\Inv(u)=\Inv(v),\Nsp(u)=\Nsp(v)$ and $\Nega(v)=\Nega(u)\cup\{u^{-1}(1)\}.$
If $i\in[n-1]$, then from $u\prec v$ and $v=s_iu$ we see that $u^{-1}(i)<u^{-1}(i+1)$, so
$v$ is the signed permutation obtained from $u$ by replacing $i$ (respectively, $\bar{i},i+1,\bar{i+1}$) by $i+1$ (respectively, $\bar{i+1},i,\bar{i}$), if they exist. So the proof follows.
\end{proof}

\vspace{3mm}
\begin{center}
\begin{tikzpicture}
\node (bar{1}bar{2}bar{3}) at (0,6.75) {{\tiny$\bar{1}\thinspace\bar{2}\thinspace \bar{3}$}};
\node (bar{1}bar{3}bar{2}) at (-1.5,5.25) {{\tiny$\bar{1}\thinspace\bar{3}\thinspace \bar{2}$}};
\node (1bar{2}bar{3}) at (0,5.25) {{\tiny$1\thinspace\bar{2}\thinspace \bar{3}$}};
\node (bar{2}bar{1}bar{3}) at (1.5,5.25) {{\tiny$\bar{2}\thinspace\bar{1}\thinspace \bar{3}$}};
\node (bar{2}bar{3}bar{1}) at (-3,3.75) {{\tiny$\bar{2}\thinspace\bar{3}\thinspace \bar{1}$}};
\node (1bar{3}bar{2}) at (-1.5,3.75) {{\tiny$1\thinspace\bar{3}\thinspace \bar{2}$}};
\node (bar{3}bar{1}bar{2}) at (0,3.75) {{\tiny$\bar{3}\thinspace\bar{1}\thinspace \bar{2}$}};
\node (2bar{1}bar{3}) at (1.5,3.75) {{\tiny$2\thinspace\bar{1}\thinspace \bar{3}$}};
\node (bar{2}1bar{3}) at (3,3.75) {{\tiny$\bar{2}\thinspace1\thinspace \bar{3}$}};
\node (bar{2}bar{3}1) at (-4.5,2.25) {{\tiny$\bar{2}\thinspace\bar{3}\thinspace1$}};
\node (2bar{3}bar{1}) at (-3,2.25) {{\tiny$2\thinspace\bar{3}\thinspace\bar{1}$}};
\node (bar{3}bar{2}bar{1}) at (-1.5,2.25) {{\tiny$\bar{3}\thinspace\bar{2}\thinspace\bar{1}$}};
\node (3bar{1}bar{2}) at (0,2.25) {{\tiny$3\thinspace\bar{1}\thinspace\bar{2}$}};
\node (bar{3}1bar{2}) at (1.5,2.25) {{\tiny$\bar{3}\thinspace1\thinspace\bar{2}$}};
\node (21bar{3}) at (3,2.25) {{\tiny$21\thinspace\bar{3}$}};
\node (bar{1}2bar{3}) at (4.5,2.25) {{\tiny$\bar{1}\thinspace2\thinspace\bar{3}$}};
\node (bar{1}bar{3}2) at (-5.25,0.75) {{\tiny$\bar{1}\thinspace\bar{3}\thinspace2$}};
\node (2bar{3}1) at (-3.75,0.75) {{\tiny$2\thinspace\bar{3}\thinspace1$}};
\node (bar{3}bar{2}1) at (-2.25,0.75) {{\tiny$\bar{3}\thinspace\bar{2}\thinspace1$}};
\node (3bar{2}bar{1}) at (-0.75,0.75) {{\tiny$3\thinspace\bar{2}\thinspace\bar{1}$}};
\node (bar{3}2bar{1}) at (0.75,0.75) {{\tiny$\bar{3}\thinspace2\thinspace\bar{1}$}};
\node (31bar{2}) at (2.25,0.75) {{\tiny$3\thinspace1\thinspace\bar{2}$}};
\node (bar{1}3bar{2}) at (3.75,0.75) {{\tiny$\bar{1}\thinspace3\thinspace\bar{2}$}};
\node (12bar{3}) at (5.25,0.75) {{\tiny$1\thinspace2\thinspace\bar{3}$}};
\node (bar{1}bar{2}3) at (-5.25,-0.75) {{\tiny$\bar{1}\thinspace\bar{2}\thinspace3$}};
\node (1bar{3}2) at (-3.75,-0.75) {{\tiny$1\thinspace\bar{3}\thinspace2$}};
\node (bar{3}bar{1}2) at (-2.25,-0.75) {{\tiny$\bar{3}\thinspace\bar{1}\thinspace2$}};
\node (3bar{2}1) at (-0.75,-0.75) {{\tiny$3\thinspace\bar{2}\thinspace1$}};
\node (bar{3}21) at (0.75,-0.75) {{\tiny$\bar{3}\thinspace2\thinspace1$}};
\node (32bar{1}) at (2.25,-0.75) {{\tiny$3\thinspace2\thinspace\bar{1}$}};
\node (bar{2}3bar{1}) at (3.75,-0.75) {{\tiny$\bar{2}\thinspace3\thinspace\bar{1}$}};
\node (13bar{2}) at (5.25,-0.75) {{\tiny$1\thinspace3\thinspace\bar{2}$}};
\node (1bar{2}3) at (-4.5,-2.25) {{\tiny$1\thinspace\bar{2}\thinspace3$}};
\node (bar{2}bar{1}3) at (-3,-2.25) {{\tiny$\bar{2}\thinspace\bar{1}\thinspace3$}};
\node (3bar{1}2) at (-1.5,-2.25) {{\tiny$3\thinspace\bar{1}\thinspace2$}};
\node (bar{3}12) at (0,-2.25) {{\tiny$\bar{3}\thinspace1\thinspace2$}};
\node (321) at (1.5,-2.25) {{\tiny$3\thinspace2\thinspace1$}};
\node (bar{2}31) at (3,-2.25) {{\tiny$\bar{2}\thinspace3\thinspace1$}};
\node (23bar{1}) at (4.5,-2.25) {{\tiny$2\thinspace3\thinspace\bar{1}$}};
\node (2bar{1}3) at (-3,-3.75) {{\tiny$2\thinspace\bar{1}\thinspace3$}};
\node (bar{2}13) at (-1.5,-3.75) {{\tiny$\bar{2}\thinspace1\thinspace3$}};
\node (312) at (0,-3.75) {{\tiny$3\thinspace1\thinspace2$}};
\node (bar{1}32) at (1.5,-3.75) {{\tiny$\bar{1}\thinspace3\thinspace2$}};
\node (231) at (3,-3.75) {{\tiny$2\thinspace3\thinspace1$}};
\node (213) at (-1.5,-5.25) {{\tiny$2\thinspace1\thinspace3$}};
\node (bar{1}23) at (0,-5.25) {{\tiny$\bar{1}\thinspace2\thinspace3$}};
\node (132) at (1.5,-5.25) {{\tiny$1\thinspace3\thinspace2$}};
\node (123) at (0,-6.75) {{\tiny$1\thinspace2\thinspace3$}};
\draw[solid,line width=0.8pt] (bar{1}bar{2}bar{3})--(bar{1}bar{3}bar{2})--(bar{2}bar{3}bar{1})--(bar{2}bar{3}1)--(bar{1}bar{3}2)--(bar{1}bar{2}3)--(1bar{2}3)--(2bar{1}3)--(213)--(123)--(132)--(231)--
(23bar{1})--(13bar{2})--(12bar{3})--(bar{1}2bar{3})--(bar{2}1bar{3})--(bar{2}bar{1}bar{3})--(bar{1}bar{2}bar{3})--(1bar{2}bar{3})--(1bar{3}bar{2})--(2bar{3}bar{1})--(2bar{3}1)
--(1bar{3}2)--(1bar{2}3);
\draw[solid,line width=0.8pt] (1bar{2}bar{3})--(2bar{1}bar{3})--(21bar{3})--(12bar{3});
\draw[solid,line width=0.8pt] (2bar{1}bar{3})--(3bar{1}bar{2})--(3bar{2}bar{1})--(3bar{2}1)--(3bar{1}2)--(2bar{1}3);
\draw[solid,line width=0.8pt] (21bar{3})--(31bar{2})--(32bar{1})--(321)--(312)--(213);
\draw[solid,line width=0.8pt] (31bar{2})--(3bar{1}bar{2});
\draw[solid,line width=0.8pt] (32bar{1})--(23bar{1});
\draw[solid,line width=0.8pt] (321)--(231);
\draw[solid,line width=0.8pt] (1bar{3}2)--(bar{1}bar{3}2);
\draw[solid,line width=0.8pt] (1bar{3}bar{2})--(bar{1}bar{3}bar{2});
\draw[solid,line width=0.8pt] (3bar{2}bar{1})--(2bar{3}bar{1});
\draw[solid,line width=0.8pt] (3bar{2}1)--(2bar{3}1);
\draw[solid,line width=0.8pt] (3bar{1}2)--(312);
\draw[dashed,line width=0.8pt] (bar{2}bar{1}bar{3})--(bar{3}bar{1}bar{2})--(bar{3}bar{2}bar{1})--(bar{3}bar{2}1)--(bar{3}bar{1}2)--(bar{2}bar{1}3)--(bar{1}bar{2}3);
\draw[dashed,line width=0.8pt] (bar{2}bar{1}3)--(bar{2}13)--(bar{1}23)--(123);
\draw[dashed,line width=0.8pt] (bar{1}23)--(bar{1}32)--(bar{2}31)--(bar{2}3bar{1})--(bar{1}3bar{2})--(13bar{2});
\draw[dashed,line width=0.8pt] (bar{1}2bar{3})--(bar{1}3bar{2});
\draw[dashed,line width=0.8pt] (bar{2}1bar{3})--(bar{3}1bar{2})--(bar{3}2bar{1})--(bar{3}21)--(bar{3}12)--(bar{2}13);
\draw[dashed,line width=0.8pt] (bar{1}32)--(132);
\draw[dashed,line width=0.8pt] (bar{2}31)--(bar{3}21);
\draw[dashed,line width=0.8pt] (bar{2}3bar{1})--(bar{3}2bar{1});
\draw[dashed,line width=0.8pt] (bar{3}12)--(bar{3}bar{1}2);
\draw[dashed,line width=0.8pt] (bar{3}1bar{2})--(bar{3}bar{1}bar{2});
\draw[dashed,line width=0.8pt] (bar{3}bar{2}1)--(bar{2}bar{3}1);
\draw[dashed,line width=0.8pt] (bar{3}bar{2}bar{1})--(bar{2}bar{3}bar{1});
\node(Figure 1.)   at (0,-7.45) {Figure $1.$ The weak order on $\mathfrak{B}_3$.};
\end{tikzpicture}
\end{center}

\begin{lemma}\label{lem:nspintr1}
Let $w$ be a signed permutation of $\mathfrak{B}_n$. If $(i,j)\in\Nsp(w)$, then $(i,\bar{j})(\bar{i},j)\in T_R(w)$.
\end{lemma}
\begin{proof}
Choose a reduced word $w=s_{i_1}s_{i_2}\cdots s_{i_k}$. Assume that $(i,j)\in\Nsp(w)$ and let $t=(i,\bar{j})(\bar{i},j)$.
Then $$s_{i_k}\prec s_{i_{k-1}}s_{i_k}\prec\cdots\prec  s_{i_1}\cdots s_{i_{k-1}}s_{i_k}=w.$$
It follows from Lemma \ref{lem:leqinn} that
$$
\emptyset=\Nsp(s_{i_k})\subseteq \Nsp(s_{i_{k-1}}s_{i_k})\subseteq\cdots\subseteq\Nsp(w),
$$
so there exists $p\in[k]$ such that
$$(i,j)\in \Nsp(s_{i_p}s_{i_{p+1}}\cdots s_{i_k})\backslash \Nsp(s_{i_{p+1}}\cdots s_{i_k}).$$
Let $u=s_{i_{k}}\cdots s_{i_{p+1}}$ and $v= s_{i_k}\cdots s_{i_{p+1}}s_{i_p}$.
Then $(i,j)\in \Nsp(v^{-1})\backslash \Nsp(u^{-1})$
and $u^{-1}$ is covered by $v^{-1}$ in the weak order. Again by Lemma \ref{lem:leqinn}, either
$$
i=u(\bar{i_p})\quad \text{and}\quad j=u(i_p+1),\quad\text{where}\quad 0<u(\bar{i_p})<u(i_p+1),
$$
or
$$
i=u(i_p+1)\quad \text{and}\quad j=u(\bar{i_p}), \quad \text{where}\quad 0<u(i_p+1)<u(\bar{i_p}).
$$
In the first case,
\begin{align*}
   t=(i,\bar{j})(\bar{i},j)=\Big(\bar{u(i_p)},\bar{u(i_p+1)}\Big)\Big(u(i_p),u(i_p+1)\Big)=us_{i_p}u^{-1},
\end{align*}
and in the second case,
\begin{align*}
    t=(i,\bar{j})(\bar{i},j)=\Big(u(i_p+1),u(i_p)\Big)\Big(\bar{u(i_p+1)},\bar{u(i_p)}\Big)=us_{i_p}u^{-1}.
\end{align*}
Thus, in both cases, we have
$$t=us_{i_p}u^{-1}=s_{i_{k}}\cdots s_{i_{p+1}}s_{i_p}s_{i_{p+1}}\cdots s_{i_{k}},$$
and hence
$$
wt=s_{i_1}\cdots s_{i_{p-1}}s_{i_{p+1}}\cdots s_{i_k}.
$$
Consequently, $\ell(wt)<\ell(w)$, so that $t\in T_R(w)$, completing the proof.
\end{proof}

\begin{lemma}[\cite{BB05}, Proposition 8.1.5]\label{lem:bbbookprop8.1.5}
The set of reflections of $\mathfrak{B}_n$ is the disjoint union
$$
T=\{(i,j)(\bar{i},\bar{j})\,|\, 1\leq i<|j|\leq n\}\biguplus\{(i,\bar{i})\,|\, i\in[n]\}.
$$
\end{lemma}

The next result provides an equivalent description of the set of right associated reflections to a signed permutation.

\begin{lemma}\label{lem:TRinvnegnsp}
Let $w$ be a signed permutation of $\mathfrak{B}_n$. Then we have a disjoint union
$$T_R(w)=\{(i,j)(\bar{i},\bar{j})\,|\,(i,j)\in\Inv(w)\}\biguplus\{(i,\bar{i})\,|\,i\in \Nega(w)\}\biguplus\{(i,\bar{j})(\bar{i},j)\,|\,(i,j)\in\Nsp(w)\}.$$
\end{lemma}
\begin{proof}
Let $A$ be the set in the right-hand side of the desired identity. By Lemma \ref{lem:bbbookprop8.1.5}, $A$ is a subset of $T$. We now show that $A\subseteq T_R(w)$.
Let $t$ be an arbitrary element of $A$. Then it suffices to show that $\ell(wt)<\ell(w)$.

If $t=(i,j)(\bar{i},\bar{j})$ for some $(i,j)\in\Inv(w)$, then $1\leq i<j\leq n$, $w_i>w_j$, and the one-line notation of $wt$ is obtained from that of $w$ by exchanging the positions of $w_i$ and $w_j$. So
$$\nega(wt)+\nsp(wt)=\nega(w)+\nsp(w)$$ and
$$
\inv(wt)=\inv(w)-2\#\{k\in [i+1,j-1]\,|\,w_j<w_k<w_i\}-1.
$$
Hence $\ell(wt)<\ell(w)$ by Eq. \eqref{eq:ell(w)b}.
If $t=(i,\bar{i})$ for some $i\in\Nega(w)$, then $w_i<0$, $(wt)(i)=\bar{w_i}$ and $wt(j)=w_j$ for all $j\in [n]$ with $j\neq i$.
Thus,
$$\nega(wt)+\nsp(wt)=\nega(w)+\nsp(w)+w_i$$
and
\begin{align*}
\inv(wt)=~&\inv(w)-\#\{k\in [i-1]\,|\,w_i<w_k<\bar{w_i}\}+\#\{k\in [i+1,n]\,|\, w_i<w_k<\bar{w_i}\}\\
\leq~&\inv(w)+\#\{k\in [i+1,n]\,|\, w_i<w_k<\bar{w_i}\}\\
\leq~&\inv(w)+\bar{w_i}-1.
\end{align*}
Again by Eq. \eqref{eq:ell(w)b}, we see that $\ell(wt)<\ell(w)$. It now follows from Lemma \ref{lem:nspintr1} that $\ell(wt)<\ell(w)$ for all $t\in A$, and hence $A\subseteq T_R(w)$.

Our task now is to show that $A=T_R(w)$. It is straightforward to verify that $A$ is a disjoint union of the three sets in the right-hand side of the desired identity.
Thus,
$$
\#A=\inv(w)+\nega(w)+\nsp(w)=\ell(w)=\#T_R(w),
$$
so that $A=T_R(w)$, and the proof follows.
\end{proof}

We are now ready to give the proof of Theorem \ref{thm:leqliffinvnegnsp}.

\begin{proof}[Proof of Theorem \ref{thm:leqliffinvnegnsp}]
If $\Inv(u)\subseteq \Inv(v)$, $\Nega(u)\subseteq\Nega(v)$ and $\Nsp(u)\subseteq\Nsp(v)$, then, by Lemma \ref{lem:TRinvnegnsp},
$T_R(u)\subseteq T_R(v)$, which together with Lemma \ref{lem:uleqlvequiv} yields that $u\leq  v$.
Conversely, assume that $u\leq v$.
Without loss of generality, we assume that $u$ is covered by $v$, that is, $\ell(v)=\ell(u)+1$ and there exists $i\in [0,n-1]$ such that
$v=s_iu$. Then the proof follows from Lemma \ref{lem:leqinn}.
\end{proof}

\section{Convex embeddings on Coxeter groups}\label{sec:opmonC}
As special order-preserving maps, parabolic convex embeddings on the symmetric group admit nice combinatorial properties \cite[Proposition 2.10]{AS05}, which were used by Aguiar and Sottile to establish the product formula of two monomial basis elements of the Malvenuto-Reutenauer Hopf algebra of permutations.
In this section, we first generalize the result of Aguiar and Sottile to any Coxeter group, and then study the convexities of a kind of nonparabolic embedding on the hyperoctahedral group. In particular, we classify the order structures of the shifted products of hyperoctahedral groups.

Throughout the rest of the paper, unless otherwise specified, $n$, $p$ and $q$ will be nonnegative integers.

\subsection{Parabolic convex embeddings on Coxeter groups}\label{subsec:opmonC}

\begin{defn}\label{defn:convexembedding}
Let $P$ and $Q$ be posets and $f:P\rightarrow Q$ a map. Then $f$ is a \emph{convex embedding} provided the following conditions are fulfilled:
\begin{enumerate}
  \item\label{defncea} $f$ is injective;
  \item\label{defnceb} for all $a,b\in P$, we have $a\leq_P b \Leftrightarrow f(a)\leq_Q f(b)$;
  \item\label{defncec} $f$ is \emph{convex}: if $f(a)\leq_Qx\leq_Q f(c)$ for some $a,c\in P$ and $x\in Q$, then there exists $b\in P$ such that $x=f(b)$.
\end{enumerate}
\end{defn}
It is remarkable that a convex embedding is always a lattice homomorphism if $P,Q$ are lattices.

\begin{lemma}\label{lem:convexlat}
Let  $P$ and $Q$ be posets and $f:P\rightarrow Q$ a convex embedding.
\begin{enumerate}
  \item\label{item:latticehomf} If $P,Q$ are lattices, then $f$ is a lattice homomorphism.
  \item\label{item:interf} If $P$ is an interval, say $p=[0_P,1_P]$, then  $f(P)=[f(0_P),f(1_P)]$.
\end{enumerate}
\end{lemma}
\begin{proof}
\eqref{item:latticehomf} Assume that $f:P\rightarrow Q$ is a convex embedding.
Let $a,b$ be any elements of $P$ and let $x=f(a)\wedge f(b)$. Then from $a\wedge b\leq a$ and $a\wedge b\leq b$ we see that
$$f(a\wedge b)\leq f(a)\quad \text{and}\quad f(a\wedge b)\leq f(b)$$
 since $f$ is order-preserving, and hence  $f(a\wedge b)\leq x\leq f(a)$.
It follows from the convexity of $f$ that there exists $c\in P$ such that $x=f(c)$, and hence $f(a\wedge b)\leq f(c) \leq f(a)$.
According to Definition \ref{defn:convexembedding}\eqref{defnceb}, we have $a\wedge b\leq c\leq a$. Similarly, we have $a\wedge b\leq c\leq b$ so that $a\wedge b= c$, which yields that \
$f(a\wedge b)=f(a)\wedge f(b)$. Simple symmetry arguments show that $f$ preserves joins.

\eqref{item:interf} Since $f$ is an order-preserving injection,
we have $f(P)\subseteq [f(0_P),f(1_P)]$. Take any $y\in Q$ with $f(0_P)\leq y\leq f(1_P)$. Then, by the convexity of $f$, there exists
$d\in P$ such that $y=f(d)\in f(P)$. Therefore,  $f(P)=[f(0_P),f(1_P)]$.
This completes the proof.
\end{proof}

Let  $J=S^A\backslash\{s_{p_1},s_{p_2},\ldots,s_{p_k}\}$ where $1\leq p_1<p_2<\cdots< p_k\leq n-1$. Fix $\xi\in \Sh(p_1,p_2-p_1,\ldots,n-p_k)$.
 Aguiar and Sottile \cite[Proposition 2.10]{AS05} showed that the map
\begin{align*}
    \rho_\xi:\mathfrak{S}_J\rightarrow \mathfrak{S}_n,\quad u \mapsto u\xi^{-1}
\end{align*}
is a convex embedding.
This is indeed true for any Coxeter groups.

\begin{theorem}\label{thm:wjembconv}
Let $(W,S)$ be any Coxeter group and $J$ a subset of $S$.
Given $\xi\in W$, let $\rho_\xi: W_J\rightarrow W$ be the map defined by $\rho_{\xi}(u)=u\xi^{-1}$. Then $\rho_\xi$ is a convex embedding if and only if $\xi\in W^J$.
\end{theorem}
\begin{proof}
Assume that $\xi\in W^J$.
The injectivity of $\rho_\xi$ is trivial since $W$ is a group.
Note that $(W^J)^{-1}:=\{w^{-1}\,|\,w\in W^J\}$ is the complete set of the minimal right coset representatives of $W_J$ in $W$.
For any $u,v$ in $W_J$, it follows from
$\ell(u\xi^{-1})=\ell(u)+\ell(\xi^{-1})=\ell(u)+\ell(\xi)$ and $\ell(v\xi^{-1})=\ell(v)+\ell(\xi)$ that
\begin{align*}
u\leq v &\Leftrightarrow \ell(v)=\ell(u)+\ell(vu^{-1})\\
&\Leftrightarrow \ell(v\xi^{-1})=\ell(u\xi^{-1})+\ell(vu^{-1})\\
&\Leftrightarrow \ell(v\xi^{-1})=\ell(u\xi^{-1})+\ell(v\xi^{-1}(u\xi^{-1})^{-1})\\
&\Leftrightarrow u\xi^{-1}\leq v\xi^{-1}\\
&\Leftrightarrow \rho_\xi(u)\leq \rho_\xi(v).
\end{align*}
It remains to show that $\rho_\xi$ is convex.  Let $u,v\in W_J$ and $w\in W$ such that $\rho_\xi (u)\leq w\leq \rho_\xi (v)$.
Since $W=W_J\cdot (W^J)^{-1}$, there exist $z\in W_J$ and $\eta\in W^J$ such that $w=z\eta^{-1}$, and hence
$$
u\xi^{-1}\leq  z\eta^{-1}\leq v\xi^{-1}.
$$
Let $\ell(v\xi^{-1})-\ell(z\eta^{-1})=k$ and $\ell(z\eta^{-1})-\ell(u\xi^{-1})=m-k$ for some nonnegative integers $k$ and $m$.
Then, by Lemma \ref{lem:uleqlvequiv}, there exists a reduced word $s_1\cdots s_ks_{k+1}\cdots s_m$ such that
\begin{align}\label{eq:azxw0J-1}
v\xi^{-1}=s_1\cdots s_ks_{k+1}\cdots s_m u\xi^{-1}\quad \ \text{and}\ \quad z\eta^{-1}=s_{k+1}\cdots s_m u\xi^{-1},
\end{align}
and hence
$$
s_1\cdots s_ks_{k+1}\cdots s_m=vu^{-1}\in W_J.
$$
Since any element of a Coxeter group has the same set of letters appearing in its reduced words, we get $s_i\in J$ for all $i\in[m]$.
Consequently, it follows from Eq.\,\eqref{eq:azxw0J-1} that
$$
\eta^{-1}\xi=z^{-1}s_{k+1}\cdots s_m u \in W_J,
$$
which together with $\eta\in W^J$ and $\xi\in W^J$ yields that $\xi=\eta$, so that $w=\rho_\xi(z)$. Thus, $\rho_u$ is convex.

Conversely, suppose that $\rho_\xi$ is a convex embedding. In particular, $\rho_\xi$ is order-preserving.  Then for any $s\in J$, it follows from $e\leq  s$ that $\xi^{-1}=\rho_\xi(e)\leq \rho_\xi(s)=s\xi^{-1}$, and hence
$$
\ell(\xi s)=\ell(s\xi^{-1})\geq \ell(\xi^{-1})=\ell(\xi)
$$
which implies that $\xi\in W^J$, completing the proof.
\end{proof}

\subsection{Nonparabolic convex embeddings on hyperoctahedral groups}
Given a sequence of positive integers $p_1,p_2,\ldots,p_k$, according to Theorem \ref{thm:wjembconv},
the embedding of the parabolic subgroup $\mathfrak{B}_{p_1}\times \mathfrak{S}_{p_2}\cdots\times \mathfrak{S}_{p_k}$ into $\mathfrak{B}_{p_1+P_2+\cdots+p_k}$
defined by
\begin{align}\label{eq:rhoxiparbloiccens}
    \rho_\xi: \mathfrak{B}_{p_1}\times \mathfrak{S}_{p_2}\cdots\times \mathfrak{S}_{p_k}\rightarrow \mathfrak{B}_{p_1+p_2+\cdots+p_k},\quad u\mapsto u\xi^{-1}
\end{align}
is convex. Here $\xi$ is any element of the set of minimal left coset representatives of $\mathfrak{B}_{p_1}\times \mathfrak{S}_{p_2}\cdots\times \mathfrak{S}_{p_k}$ in
$\mathfrak{B}_{p_1+P_2+\cdots+p_k}$.

There exists another similar construction for hyperoctahedral groups. Given positive integers $p$ and $q$,
consider the shifted product $\mathfrak{B}_{p}\times\mathfrak{B}_{q}$ to be the subgroup of $\mathfrak{B}_{p+q}$, where we identify the element $(u,v)\in\mathfrak{B}_{p}\times\mathfrak{B}_{q}$
with the signed permutation $u\times v$ of $\mathfrak{B}_{p+q}$ defined by
\begin{align*}
    u\times v:=(u_1,\ldots,u_p, v_1[p],\ldots,v_{q}[p]).
\end{align*}
Here for nonzero integer $a$ we use the notation
\begin{align*}
a[p]:=a+{\rm sgn}(a)\cdot p=\begin{cases}
a+p,&a>0,\\
a-p,&a<0.
\end{cases}
\end{align*}
For example, $2\bar{4}\,\bar{3}1[3]=5\bar{7}\,\bar{6}4$ and $2\bar{1}3\times 2\bar{4}\,\bar{3}1=2\bar{1}35\bar{7}\,\bar{6}4$.
The set of $(p,q)$-shuffles
$$
\Sh(p,q)=\{\xi\in \mathfrak{S}_{p+q}\,|\,\xi_1<\cdots<\xi_{p},\,\xi_{p+1}<\cdots<\xi_{p+q}\}
$$
is a set of left coset representatives for $\mathfrak{B}_p\times\mathfrak{B}_{q}$ in the group $\mathfrak{B}_{p+q}$.
Thus we have a bijection
\begin{align*}
\sigma: \Sh(p,q)\times \mathfrak{B}_{p}\times \mathfrak{B}_{q}\rightarrow\mathfrak{B}_{p+q},\quad (\xi,u,v)\mapsto \xi(u\times v).
\end{align*}
More generally, given positive integers $p_1,p_2,\ldots,p_k$, the shuffle set $\Sh(p_1,p_2,\ldots,p_k)$ is a set of the left coset representatives for the subgroup
$\mathfrak{B}_{p_1}\times \mathfrak{B}_{p_2}\times \cdots\times \mathfrak{B}_{p_k}$
in $\mathfrak{B}_{p_1+p_2+\cdots +p_k}$.

We now study the order structure of the subgroup $\mathfrak{B}_{p_1}\times \mathfrak{B}_{p_2}\times \cdots\times \mathfrak{B}_{p_k}$.
To mimic the map given by \eqref{eq:rhoxiparbloiccens}, for any $\xi\in \Sh(p_1,p_2,\ldots,p_k)$, we define the nonparabolic embedding
\begin{align*}
\tau_{\xi}:\mathfrak{B}_{p_1}\times\mathfrak{B}_{p_2}\times \cdots\times \mathfrak{B}_{p_k}\rightarrow \mathfrak{B}_{p_1+p_2+\cdots+p_k},\quad u\mapsto u\xi^{-1}.
\end{align*}
Given a subset $K$ of $[n]$, let
$$
\mathfrak{B}_{n,K}:=\{w\in \mathfrak{B}_n\,|\, \Nega(w)=K\}.
$$
Let $L_i\subseteq [p_i]$ where $i=2,3\ldots,k$,
the set $\mathfrak{B}_{p_1}\times  \mathfrak{B}_{p_2,L_2}\times\cdots\times  \mathfrak{B}_{p_k,L_k}$ is called the \emph{$(L_2,L_3,\ldots,L_k)$-component} of
$\mathfrak{B}_{p_1+p_2+\cdots+p_k}$.
We show below by induction on $k$ that the restriction of $\tau_{\xi}$ to each component of $\mathfrak{B}_{p_1}\times\mathfrak{B}_{p_2}\times \cdots\times \mathfrak{B}_{p_k}$ is a convex embedding.

The following statement is straightforward and the proof is omitted.
\begin{lemma}\label{lem:wxistbw12}
Let $w\in\mathfrak{B}_n$, and let $\xi=w\cdot (\sts(w_1\cdots w_p)\times\sts(w_{p+1}\cdots w_{n}))^{-1}$. Then $\xi\in \Sh(p,n-p)$ and
$\sigma^{-1}(w)=(\xi,\sts(w_1\cdots w_p),\sts(w_{p+1}\cdots w_{n}))$.
\end{lemma}

There is a simple combinatorial rule to get the shuffle $\xi=w\cdot (\sts(w_1\cdots w_p)\times\sts(w_{p+1}\cdots w_{n}))^{-1}$.
Namely, $\xi$ is obtained from $w$
by first rearranging  in increasing order the absolute values $|w_1|,|w_2|$, $\cdots$, $|w_{p}|$ so that they appear in the places $1,\ldots,p$, and then rearranging the absolute values $|w_{p+1}|,|w_{p+2}|,\ldots,|w_{n}|$ so that they appear in increasing order in the places $p+1,p+2,\ldots,n$.
For example,
$$
4\bar{2}\,\bar{6}17\bar{3}5\bar{8}=12463578\cdot(3\bar{2}\,\bar{4}1\times 3\bar{1}2\bar{4}).
$$
\begin{lemma}\label{lem:stsuv}
Let $u\leq v$ in $\mathfrak{B}_n$. Then for any positive integers $i$ and $j$ with $1\leq i<j\leq n$, we have $\sts(u_{i}u_{i+1}\cdots u_{j})\leq \sts(v_{i}v_{i+1}\cdots{v_j})$.
\end{lemma}
\begin{proof}
Let $w\in \mathfrak{B}_n$ and write $w_{ij}=w_iw_{i+1}\cdots w_{j}$.
Observe that
\begin{align*}
    \Inv(\sts(w_{ij}))=&(\Inv(w)\cap[i,j]\times [i,j])-(i-1,i-1),\\
    \Nega(\sts(w_{ij}))=&(\Nega(w)\cap[i,j])-i+1
    \intertext{and}
    \Nsp(\sts(w_{ij}))=&(\Nsp(w)\cap[i,j]\times [i,j])-(i-1,i-1).
\end{align*}
Thus, if $u\leq v$, then $\sts(u_{ij})\leq \sts(v_{ij})$ by Theorem \ref{thm:leqliffinvnegnsp}.
\end{proof}

Given a signed permutation $w\in\mathfrak{B}_n$, let
\begin{align*}
\overline{\Nega}(w):=\{i\in[n]\,|\,i\not\in\Nega(w)\},\quad
    \widetilde{\Inv}(w)):=\{(j,i)\in[n]\times[n]\,|\,(i,j)\in \Inv(w)\}
\end{align*}
and
\begin{align*}
    \overline{\Inv}(w):=\{(i,j)\in[n]\times[n]\,|\,i<j,\ (i,j)\not\in \Inv(w)\}.
\end{align*}
Then we have the following key observation about the elements belonging to the image of the nonparabolic embedding $\tau_{\xi}$.

\begin{lemma}\label{lem:uvleqw}
Let $w=(u\times v) \xi^{-1}$ where $\xi\in \Sh(p,q)$, $u\in\mathfrak{B}_p$ and $v\in\mathfrak{B}_q$. Then we have the following disjoint unions
\begin{align}
\Nega(w)={}&\xi\cdot \Nega(u)\biguplus\xi\cdot\Big(p+\Nega(v)\Big),\label{eq:item-nega(w)}\\
\Inv(w)={}&\xi\cdot \Inv(u) \biguplus\xi\cdot\Big((p,p)+\Inv(v)\Big)\biguplus\xi\cdot \left([p]\times (p+\Nega(v))\cap \overline{\Inv}(\xi)\right)\notag\\
&\biguplus\xi\cdot\left((p+\overline{\Nega}(v))\times [p]\cap\widetilde{\Inv}(\xi)\right)\label{eq:item-Inv(w)}
\intertext{and}
\Nsp(w)={}&\xi\cdot \Nsp(u) \biguplus\xi\cdot\Big((p,p)+\Nsp(v)\Big)\notag
\biguplus\xi\cdot\left([p]\times (p+\Nega(v))\cap\overline{\Inv}(\xi)\right)\\
&\biguplus\xi\cdot\left((p+\Nega(v))\times[p]\cap\widetilde{\Inv}(\xi)\right)\label{eq:item-Nsp(w)}.
\end{align}
\end{lemma}
\begin{proof}
Let $i$ be a positive integer and assume that $i=\xi_{r}$ for some $r\in [p+q]$. Then
\begin{align*}
    i\in\Nega(w) \Leftrightarrow{}&  w(\xi_{r})<0\Leftrightarrow (u\times v)(r)<0\Leftrightarrow r\in \Nega(u)\cup (p+\Nega(v))\\
    \Leftrightarrow{}& i\in \xi\cdot\left(\Nega(u)\cup(p+\Nega(v))\right),
\end{align*}
so Eq.\,\eqref{eq:item-nega(w)} follows.

Now consider Eqs.\,\eqref{eq:item-Inv(w)} and \eqref{eq:item-Nsp(w)}.
Let $i$ and $j$ be positive integers with $1\leq i<j\leq p+q$, and assume that $i=\xi_r$ and $j=\xi_s$ for some $r,s\in[p+q]$.
Then $\xi_r<\xi_s$. We use $(i,j)=\xi\cdot (r,s)$ to indicate that $i=\xi_r$ and $j=\xi_s$.
It follows from $\xi\in \Sh(p,q)$ that there are four cases need to be considered.

If $r,s\in [p]$, then
\begin{align*}
   w_i=w(\xi_r)=(u\times v)(r)=u_r\quad\text{and}\quad w_j=w(\xi_s)=(u\times v)(s)=u_s.
\end{align*}
It follows from $\xi_r<\xi_s$ that $1\leq r<s\leq p$, so $(r,s)\in\overline{\Inv}(\xi)$, and hence
\begin{align*}
(i,j)\in\Inv(w)\Leftrightarrow{}w_i>w_j\Leftrightarrow{}u_r>u_s
\Leftrightarrow {}(r,s)\in \Inv(u)\cap \overline{\Inv}(\xi) \Leftrightarrow{} (i,j)\in \xi\cdot \left(\Inv(u)\cap \overline{\Inv}(\xi)\right)
\end{align*}
and
\begin{align*}
(i,j)\in\Nsp(w)&\Leftrightarrow{}w_i+w_j<0\Leftrightarrow{}u_r+u_s<0\\
&\Leftrightarrow {}(r,s)\in \Nsp(u)\cap\overline{\Inv}(\xi)
\Leftrightarrow{} (i,j)\in \xi\cdot \left(\Nsp(u)\cap \overline{\Inv}(\xi)\right).
\end{align*}
Since
$$\Inv(u)\subseteq\{(i,j)\,|\,1\leq i<j\leq p\}\subseteq\overline{\Inv}(\xi)\quad \text{and} \quad\Nsp(u)\subseteq\{(i,j)\,|\,1\leq i<j\leq p\}\subseteq\overline{\Inv}(\xi),$$
in this case we obtain that
\begin{align}\label{eq:invnspij1p}
(i,j)\in\Inv(w)\Leftrightarrow{} (i,j)\in \xi\cdot \Inv(u)\quad\text{and}\quad
(i,j)\in\Nsp(w)\Leftrightarrow{} (i,j)\in \xi\cdot \Nsp(u).
\end{align}

If $r,s\in[p+1,p+q]$, then
\begin{align*}
   w_i=w(\xi_r)=(u\times v)(r)=v_{r-p}[p]\quad\text{and}\quad w_j=w(\xi_s)=(u\times v)(s)=v_{s-p}[p].
\end{align*}
From $\xi_r<\xi_s$ we see that $p+1\leq r<s\leq p+q$ and $(r,s)\in\overline{\Inv}(\xi)$,  so that
\begin{align*}
(i,j)\in\Inv(w)&\Leftrightarrow{}w_i>w_j\Leftrightarrow{}v_{r-p}[p]>v_{s-p}[p]
\Leftrightarrow {}(r-p,s-p)\in \Inv(v)\\
&\Leftrightarrow (r,s)\in \left((p,p)+\Inv(v)\right)\cap \overline{\Inv}(\xi)\\
&\Leftrightarrow (i,j)\in \xi\cdot \left(((p,p)+\Inv(v))\cap \overline{\Inv}(\xi)\right)
\end{align*}
and
\begin{align*}
(i,j)\in\Nsp(w)&\Leftrightarrow{}w_i+w_j<0\Leftrightarrow{}v_{r-p}[p]+v_{s-p}[p]<0\Leftrightarrow{}v_{r-p}+v_{s-p}<0\\
&\Leftrightarrow {}((r,s)\in \Big((p,p)+\Nsp(v)\Big)\cap \overline{\Inv}(\xi)\\
&\Leftrightarrow (i,j)\in \xi\cdot \left(((p,p)+\Nsp(v))\cap \overline{\Inv}(\xi)\right).
\end{align*}
In this case we have
\begin{align*}(p,p)+\Inv(v)&\subseteq\{(i,j)\,|\,p+1\leq i<j\leq p+q\}\subseteq \overline{\Inv}(\xi)
\intertext{and}
(p,p)+\Nsp(v)&\subseteq\{(i,j)\,|\,p+1\leq i<j\leq p+q\}\subseteq \overline{\Inv}(\xi),
\end{align*}
so that
\begin{align}\label{eq:invij1pqp}
(i,j)\in\Inv(w) \Leftrightarrow (i,j)\in \xi\cdot \Big((p,p)+\Inv(v)\Big)
\end{align}
and
\begin{align}\label{eq:nspij1ppq}
(i,j)\in\Nsp(w)\Leftrightarrow (i,j)\in \xi\cdot \Big((p,p)+\Nsp(v)\Big).
\end{align}

If $r\in[p]$ and $s\in[p+1,p+q]$, then $(r,s)\in \overline{\Inv}(\xi)$ since $\xi_r<\xi_s$, and we have
$$w_i=(u\times v)(r)=u_{r}\quad \text{and}\quad w_j=(u\times v)(s)=v_{s-p}[p].$$
So
\begin{align}\label{eq:invpp+1p+q}
    (i,j)\in\Inv(w)&\Leftrightarrow{}w_i>w_j\Leftrightarrow{}u_{r}>v_{s-p}[p]\Leftrightarrow{}r\in[p],\ v_{s-p}<0\notag\\
    &\Leftrightarrow{}(r,s)\in [p]\times (p+\Nega(v))\cap\overline{\Inv}(\xi)\notag\\
    &\Leftrightarrow{}(i,j)\in \xi\cdot\left([p]\times (p+\Nega(v))\cap\overline{\Inv}(\xi)\right)
\end{align}
and
\begin{align}\label{eq:negapp+1p+q}
    (i,j)\in\Nsp(w)&\Leftrightarrow{}w_i+w_j<0\Leftrightarrow{}u_{r}+v_{s-p}[p]<0\Leftrightarrow{}r\in[p],\ v_{s-p}<0\notag\\
    &\Leftrightarrow{}(r,s)\in [p]\times (p+\Nega(v))\cap\overline{\Inv}(\xi)\notag\\
    &\Leftrightarrow{}(i,j)\in \xi\cdot\left([p]\times (p+\Nega(v))\cap\overline{\Inv}(\xi)\right).
\end{align}

If $r\in[p+1,p+q]$ and $s\in[p]$, then $(r,s)\in \widetilde{\Inv}(\xi)$,
$w_i=v_{r-p}[p]$ and $w_j=u_{s}$. So
\begin{align}\label{eq:invpp+1p+qp}
    (i,j)\in\Inv(w)&\Leftrightarrow{}w_i>w_j\Leftrightarrow{}v_{r-p}[p]>u_{s}\Leftrightarrow{}s\in[p],\ v_{r-p}>0\notag\\
    &\Leftrightarrow{}(r,s)\in (p+\overline{\Nega}(v))\times[p]\cap \widetilde{\Inv}(\xi)\notag\\
    &\Leftrightarrow{}(i,j)\in \xi\cdot\left((p+\overline{\Nega}(v))\times[p])\cap \widetilde{\Inv}(\xi)\right)
\end{align}
and
\begin{align}\label{eq:negapp+1p+qp}
    (i,j)\in\Nsp(w)&\Leftrightarrow{}w_i+w_j<0\Leftrightarrow{}v_{r-p}[p]+u_{s}<0\Leftrightarrow{}s\in[p],\ v_{r-p}<0\notag\\
    &\Leftrightarrow{}(r,s)\in (p+\Nega(v))\times[p]\cap \widetilde{\Inv}(\xi)\notag\\
    &\Leftrightarrow{}(i,j)\in \xi\cdot\left((p+\Nega(v))\times[p] \cap \widetilde{\Inv}(\xi)\right).
\end{align}

Eq.\,\eqref{eq:item-Inv(w)} now follows from Eqs.\,\eqref{eq:invnspij1p}, \eqref{eq:invij1pqp}, \eqref{eq:invpp+1p+q} and \eqref{eq:invpp+1p+qp}, while
Eq.\,\eqref{eq:item-Nsp(w)} follows from Eqs.\,\eqref{eq:invnspij1p}, \eqref{eq:nspij1ppq}, \eqref{eq:negapp+1p+q} and \eqref{eq:negapp+1p+qp}.
\end{proof}

\begin{coro}\label{coro:uuleqvv}
Let $u,u'\in\mathfrak{B}_p$ and $v,v'\in\mathfrak{B}_q$. Then $u\times v\leq u'\times v'$ if and only if $u\leq u'$ and $v\leq v'$.
In this case we have
$$\ell(u'\times v')-\ell(u\times v)= \ell(u')-\ell(u)+\ell(v')-\ell(v)+2p(\nega(v')-\nega(v)).$$
\end{coro}
\begin{proof}
Putting $\xi=1_{p+q}$ in Lemma \ref{lem:uvleqw} yields that
\begin{align*}
  \Nega(u\times v)=&\,\Nega(u)\biguplus(p+\Nega(v)),\\
  \Inv(u\times v)=&\,\Inv(u)\biguplus((p,p)+\Inv(v))\biguplus[p]\times \left(p+\Nega(v)\right)
  \intertext{and}
  \Nsp(u\times v)=&\,\Nsp(u)\biguplus\left((p,p)+\Nsp(v)\right)\biguplus [p]\times \left(p+\Nega(v)\right).
\end{align*}
Then the proof follows from Eq.\,\eqref{eq:ell(w)b} and Theorem \ref{thm:leqliffinvnegnsp}.
\end{proof}

The following proposition shows that
each component of the shifted product $\mathfrak{B}_p\times\mathfrak{B}_q$ can be convex embedded into  $\mathfrak{B}_{p+q}$ under the restriction of the map $\tau_{\xi}$.

\begin{prop}\label{prop:convexuvupq}
Let $\xi\in \Sh(p,q)$.
Then for any $K\subseteq [q]$, $\tau_{\xi}$ is a convex embedding when restricted to
$\mathfrak{B}_p\times  \mathfrak{B}_{q,K}$. In particular, this restriction preserves meets and joins.
\end{prop}
\begin{proof}
It is clear that $\tau_{\xi}$ is injective.
Theorem \ref{thm:leqliffinvnegnsp} and Lemma \ref{lem:uvleqw} guarantee that  $\tau_{\xi}$ is order-preserving.
It remains to prove that the restriction of $\tau_\xi$ to $\mathfrak{B}_p\times  \mathfrak{B}_{q,K}$ is convex.
To see this, set $w\in \mathfrak{B}_{p+q}$, $u',u''\in \mathfrak{B}_p$ and $v',v''\in\mathfrak{B}_{q,K}$ such that $\tau_{\xi}(u'\times v')\leq w\leq \tau_{\xi}(u''\times v'')$.
Let $w'=\tau_{\xi}(u'\times v')$ and $w''=\tau_{\xi}(u''\times v'')$. Then $w'\leq w\leq w''$.
According to Lemma \ref{lem:wxistbw12}, there exist $\eta\in \Sh(p,q)$, $u\in \mathfrak{B}_p$ and $v\in \mathfrak{B}_q$ such that
$w^{-1}=\eta(u^{-1}\times v^{-1})$, that is, $w=(u\times v)\eta^{-1}$.
We claim that $\xi=\eta$.
Suppose to the contrary that $\xi\neq\eta$. Then there exist $i,r\in[p]$
and $j,s\in[p+1,p+q]$ such that $\xi_i=\eta_s$ and $\xi_j=\eta_r$.

We first consider the case $\xi_i>\xi_j$. Then $(i,j)\in\Inv(\xi)$
and $(r,s)\in \overline{\Inv}(\eta)$. Comparing Eqs.\,\eqref{eq:item-Inv(w)} and \eqref{eq:item-Nsp(w)} yields that
\begin{align*}
    (\eta_r,\eta_s)\in \Inv(w)&\Leftrightarrow (r,s)\in [p]\times (p+\Nega(v))\cap\overline{\Inv}(\eta)\Leftrightarrow (\eta_r,\eta_s)\in \Nsp(w).
\end{align*}
Since $(\xi_j,\xi_i)=(\eta_r,\eta_s)$, we have
$$(\xi_j,\xi_i)\in \Inv(w)\Leftrightarrow (\xi_j,\xi_i)\in \Nsp(w).$$
Again by Eqs.\,\eqref{eq:item-Inv(w)} and \eqref{eq:item-Nsp(w)}, if $j\in p+\Nega(v')$, then
\begin{align*}
    (j,i)\in (p+\Nega(v'))\times [p]\cap\widetilde{\Inv}(\xi)&\Rightarrow (\xi_j,\xi_i)\in \Nsp(w')\subseteq\Nsp(w)\\
    &\Rightarrow (\xi_j,\xi_i)\in \Inv(w)\subseteq\Inv(w'')\\
    &\Rightarrow j\in p+\overline{\Nega}(v'').
\end{align*}
But $\Nega(v')=\Nega(v'')$, so $j\in p+\overline{\Nega}(v')$,
a contradiction. Similarly, if $j\in p+\overline{\Nega}(v')$, then
\begin{align*}
    (j,i)\in (p+\overline{\Nega}(v'))\times [p]\cap\widetilde{\Inv}(\xi)&\Rightarrow (\xi_j,\xi_i)\in \Inv(w')\subseteq\Inv(w)\\
    &\Rightarrow (\xi_j,\xi_i)\in \Nsp(w)\subseteq\Nsp(w'')\\
    &\Rightarrow j\in p+\Nega(v'')\\
    &\Rightarrow j\in p+\Nega(v'),
\end{align*}
a contradiction.

Thus, we must have $\xi_i<\xi_j$, and hence $(i,j)\in\overline{\Inv}(\xi)$
and $(r,s)\in \Inv(\eta)$.
According to Eqs.\,\eqref{eq:item-Inv(w)} and \eqref{eq:item-Nsp(w)}, we have
\begin{align*}
    (\xi_i,\xi_j)\in \Inv(w')&\Leftrightarrow (i,j)\in [p]\times (p+\Nega(v'))\cap\overline{\Inv}(\xi)\Leftrightarrow (\xi_i,\xi_j)\in \Nsp(w')
    \intertext{and}
    (\xi_i,\xi_j)\in \Inv(w'')&\Leftrightarrow (i,j)\in [p]\times (p+\Nega(v''))\cap\overline{\Inv}(\xi)\Leftrightarrow (\xi_i,\xi_j)\in \Nsp(w'').
\end{align*}
It then follows from $\Nega(v')=\Nega(v'')$ that
\begin{align*}
    (\xi_i,\xi_j)\in \Inv(w')\Leftrightarrow (\xi_i,\xi_j)\in \Nsp(w')\Leftrightarrow
    (\xi_i,\xi_j)\in \Inv(w'')\Leftrightarrow(\xi_i,\xi_j)\in \Nsp(w'').
\end{align*}
Consequently, in view of Eqs.\,\eqref{eq:item-Inv(w)} and \eqref{eq:item-Nsp(w)}, if $s\in p+\Nega(v)$, then
\begin{align*}
(s,r)\in(p+\Nega(v)) \times[p]\cap\widetilde{\Inv}(\eta)&\Rightarrow (\xi_i,\xi_j)=(\eta_s,\eta_r)\in \Nsp(w)\subseteq\Nsp(w'')\\
&\Rightarrow (\eta_s,\eta_r)=(\xi_i,\xi_j)\in \Inv(w')\subseteq\Inv(w)\\
&\Rightarrow s\in p+\overline{\Nega}(v);
\end{align*}
if $s\in p+\overline{\Nega}(v)$, then
\begin{align*}
(s,r)\in(p+\overline{\Nega}(v)) \times[p]\cap\widetilde{\Inv}(\eta)&\Rightarrow (\xi_i,\xi_j)=(\eta_s,\eta_r)\in \Inv(w)\subseteq\Inv(w'')\\
&\Rightarrow (\eta_s,\eta_r)=(\xi_i,\xi_j)\in \Nsp(w')\subseteq\Nsp(w)\\
&\Rightarrow s\in p+\Nega(v).
\end{align*}
Thus, in both cases, we reach a contradiction. So $\xi=\eta$.

Now we have $(u'\times v')\xi^{-1}\leq (u\times v)\xi^{-1}\leq (u''\times v'')\xi^{-1}$.
Since $\Nega(v')=\Nega(v'')=K$, we see that $\Nega(v)=K$ by Eq.\,\eqref{eq:item-nega(w)}, and hence $u\times v\in \mathfrak{B}_p\times  \mathfrak{B}_{q,K}$.
Therefore, $\tau_{\xi}$ is a convex embedding when restricted to $\mathfrak{B}_p\times  \mathfrak{B}_{q,K}$.
\end{proof}

\begin{coro}\label{coro:uesu'leqvtimv'}
Let $K\subseteq[q]$. If $u\leq u'$ in $\mathfrak{B}_p$ and $v\leq v'$ in $\mathfrak{B}_{q,K}$, then $[u,u']\times[v,v']=[u\times v,u'\times v']$ in $\mathfrak{B}_{p+q}$.
\end{coro}
\begin{proof}
Since $[x,y]\cong[1_{p+q},yx^{-1}]$ for any $x,y\in\mathfrak{B}_{p+q}$, we have $[u\times v,u'\times v']\cong [1_{p+q},u'u^{-1}\times v'v^{-1}]$.
From $v, v'\in \mathfrak{B}_{q,K}$ we see that $v'v^{-1}\in \mathfrak{B}_{q,\emptyset}$, so it suffices to assume that $u\in \mathfrak{B}_{p}$ and $v\in \mathfrak{B}_{q,\emptyset}$ and show that
$$
[1_p,u]\times[1_q,v]=[1_{p+q},u\times v].
$$
By Corollary \ref{coro:uuleqvv}, $[1_p,u]\times[1_q,v]\subseteq[1_{p+q},u\times v]$.
For the other inclusion, pick $w\in [1_{p+q},u\times v]$.
In view of Proposition \ref{prop:convexuvupq}, the map $(u,v)\mapsto u\times v$ is a convex embedding from $\mathfrak{B}_p\times  \mathfrak{B}_{q,\emptyset}$ to
$\mathfrak{B}_{p+q}$, so $w=u''\times v''$ for some $(u'',v'')\in\mathfrak{B}_p\times  \mathfrak{B}_{q,\emptyset}$.
Again by  Corollary \ref{coro:uuleqvv} we have $u''\in[1_p,u]$ and $v''\in[1_q,v']$, so that $w\in[1_p,u]\times[1_q,v]$.
Thus, $[1_p,u]\times[1_q,v]=[1_{p+q},u\times v]$, completing the proof.
\end{proof}

We now are in a position to give the main result of this subsection.
The proof relies upon the following fact given in \cite{AS05}.
Let $p_1,p_2,\ldots,p_k$ be positive integers with sum $p$. Then for any positive integer $q$, the map
\begin{align}\label{eq:shpqp1pk}
\Sh(p,q)\times \Sh(p_1,p_2,\ldots,p_k)\rightarrow \Sh(p_1,p_2,\ldots,p_k,q),\quad (\zeta,\zeta')\mapsto \zeta (\zeta'\times 1_q)
\end{align}
is a bijection.

\begin{theorem}\label{thm:convexuvupqg}
Let $p_1,p_2,\ldots,p_k$ be positive integers where $k\geq2$.
\begin{enumerate}
\item\label{lemma-item:abpqukg} We have
\begin{align*}
\mathfrak{B}_{p_1}\times\mathfrak{B}_{p_2}\times\cdots\times\mathfrak{B}_{p_k}=
\biguplus_{L_i\subseteq [p_i]\atop i=2,\ldots, k}\mathfrak{B}_{p_1}\times  \mathfrak{B}_{p_2,L_2}\times\cdots\times  \mathfrak{B}_{p_k,L_k}
\intertext{and}
\mathfrak{B}_{p_1+p_2+\cdots+p_k}=\biguplus_{\xi\in\Sh(p_1,p_2,\ldots,p_k)}\biguplus_{L_i\subseteq [p_i]\atop i=2,\ldots, k}
\left(\mathfrak{B}_{p_1}\times  \mathfrak{B}_{p_2,L_2}\times\cdots\times  \mathfrak{B}_{p_k,L_k}\right)\cdot\xi^{-1}.
\end{align*}
\item\label{lemma-item:abpqsublg} For any  $L_i\subseteq [p_i]$, $i=2,\ldots, k$,
the subset $\mathfrak{B}_{p_1}\times  \mathfrak{B}_{p_2,L_2}\times\cdots\times  \mathfrak{B}_{p_k,L_k}$
is a sublattice of $\mathfrak{B}_{p_1+p_2+\cdots+p_k}$ isomorphic to $\mathfrak{B}_{p_1}\times  \mathfrak{S}_{p_2}\times\cdots\times  \mathfrak{S}_{p_k}$.
\item\label{lemma-item:abpqccomg} Let $L_i,L_i'\subseteq [p_i]$, $i=2,\ldots, k$, where there is at least one $i$ such that $L_i\neq L_i'$.
If $w\in\mathfrak{B}_{p_1}\times  \mathfrak{B}_{p_2,L_2}\times\cdots\times  \mathfrak{B}_{p_k,L_k}$ and  $w'\in\mathfrak{B}_{p_1}\times  \mathfrak{B}_{p_2,L_2'}\times\cdots\times  \mathfrak{B}_{p_k,L_k'}$ such that $w<w'$, then $\ell(w')-\ell(w)\geq3$.
\item\label{lemma-item:tauximeetjoing} For any $\xi\in \Sh(p_1,p_2,\ldots,p_k)$ and $L_i\subseteq [p_i]$, $i=2,\ldots, k$, the map $\tau_{\xi}$ is a convex embedding when restricted to
$\mathfrak{B}_{p_1}\times  \mathfrak{B}_{p_2,L_2}\times\cdots\times  \mathfrak{B}_{p_k,L_k}$. In particular, this restriction preserves meets and joins.
\end{enumerate}
\end{theorem}
\begin{proof}
\eqref{lemma-item:abpqukg} Since $\mathfrak{B}_{p_i}=\biguplus_{L_i\subseteq [p_i]}\mathfrak{B}_{p_i,L_i}$ for $i=2,3,\ldots,k$, the first identity is true.
Note that the set $\{\xi^{-1}\,|\,\xi\in\Sh(p_1,p_2,\ldots,p_k)\}$ is a set of the right coset representatives for the subgroup $\mathfrak{B}_{p_1}\times  \mathfrak{B}_{p_2}\times\cdots\times  \mathfrak{B}_{p_k}$ in $\mathfrak{B}_{p_1+p_2+\cdots+p_k}$. Hence the second identity holds.

\eqref{lemma-item:abpqsublg} The proof is by induction on $k$.
 We first show that for any subset $L_2$ of $[p_2]$, $\mathfrak{B}_{p_1}\times  \mathfrak{B}_{p_2,L_2}$ is a sublattice of $\mathfrak{B}_{p_1+p_2}$ isomorphic to
 $\mathfrak{B}_{p_1}\times  \mathfrak{S}_{p_2}$.
Suppose that
$$L_2=\{i_1<i_2<\cdots<i_r\}\quad \text{and}\quad [p_2]\backslash L_2=\{i_{r+1}<i_{r+2}<\cdots<i_{p_2}\}.$$
Let $v',v''\in \mathfrak{B}_{p_2}$ where
\begin{align*}
    v'_{i_1}=\bar{r},\ v'_{i_2}=\bar{r-1},\ \cdots,\ v'_{i_r}=\bar{1},\ v'_{i_{r+1}}=r+1,\ v'_{i_{r+2}}=r+2,\ \cdots,\ v'_{i_{p_2}}=p_2,
    \end{align*}
and
\begin{align*}
    v''_{i_1}=\bar{p_2-r+1},\ v''_{i_2}=\bar{p_2-r+2},\ \cdots,\ v''_{i_r}=\bar{p_2},\ v''_{i_{r+1}}=p_2-r,\ v''_{i_{r+2}}=p_2-r-1,\ \cdots,\ v''_{i_{p_2}}=1.
\end{align*}
It is now straightforward to check that $\mathfrak{B}_{p_2,L_2}=[v',v'']$ by Theorem \ref{thm:leqliffinvnegnsp}.

Since $\mathfrak{B}_{p_1}=[1_{p_1},\omega_{p_1}]$ where $\omega_{p_1}$ is the largest element of $\mathfrak{B}_{p_1}$, it follows from Corollary \ref{coro:uesu'leqvtimv'} that
$\mathfrak{B}_{p_1}\times  \mathfrak{B}_{{p_2},L_2}=[1_{p_1}\times v',\omega_{p_1}\times v'']$,
which is a sublattice of $\mathfrak{B}_{p_1+p_2}$.
It is easily seen that $v''v'^{-1}=p_2(p_2-1)\cdots 1$ which is the largest element of $\mathfrak{S}_{p_2}$, so $\mathfrak{B}_{p_2,L_2}\cong[1_{p_2},v''v'^{-1}]=\mathfrak{S}_{p_2}$,
and hence
$$\mathfrak{B}_{p_1}\times  \mathfrak{B}_{p_2,L_2}\simeq[1_{p_1}\times 1_{p_2},\omega_{p_1}\times (p_2(p_2-1)\cdots 1)]=\mathfrak{B}_{p_1}\times \mathfrak{S}_{p_2}.$$

Now assume for $k-1$. Then $\mathfrak{B}_{p_k,L_k}$ and $\mathfrak{B}_{p_1}\times  \mathfrak{B}_{p_2,L_2}\times\cdots\times  \mathfrak{B}_{p_{k-1},L_{k-1}}$ are intervals,
and hence so does $\mathfrak{B}_{p_1}\times  \mathfrak{B}_{p_2,L_2}\times\cdots\times  \mathfrak{B}_{p_k,L_k}$ by Corollary \ref{coro:uesu'leqvtimv'}. It is routine to show that this interval is isomorphic to the parabolic subgroup $\mathfrak{B}_{p_1}\times  \mathfrak{S}_{p_2}\times\cdots\times  \mathfrak{S}_{p_k}$.

\eqref{lemma-item:abpqccomg} Suppose
$w=u^{(1)}\times u^{(2)}\times\cdots\times u^{(k)}$ and
$w'=u'^{(1)}\times u'^{(2)}\times\cdots\times u'^{(k)}$, where $u^{(1)},u'^{(1)}\in \mathfrak{B}_{p_1}$, $u^{(i)}\in \mathfrak{B}_{p_i,L_i}$ and $u'^{(i)}\in \mathfrak{B}_{p_i,L_i'}$ for
$i=2,\ldots,k$.
Since $w<w'$, it follows from Corollary \ref{coro:uuleqvv} that $u^{(i)}\leq u'^{(i)}$, $i=1,2,\ldots,k$. Then, by Theorem \ref{thm:leqliffinvnegnsp}, we have $L_i\subseteq L'_{i}$,
$i=2,\ldots,k$.
Choose $j$ maximal so that $L_j\neq L'_j$. Then $j\geq2$, $L_j\subsetneqq  L'_j$ and $\ell(u'^{(j)})-\ell(u^{(j)})\geq1$.
From Corollary \ref{coro:uuleqvv} it follows that
\begin{align*}
    \ell(w')-\ell(w)\geq\, & \sum_{i=1}^k\left(\ell(u'^{(i)})-\ell(u^{(i)})\right)+2(p_1+\cdots+p_{j-1})(|L'_j|-|L_j|)\\
    \geq\, & \ell(u'^{(j)})-\ell(u^{(j)})+2(p_1+\cdots+p_{j-1})\\
    \geq\, & 3.
\end{align*}

\eqref{lemma-item:tauximeetjoing} Let $\xi\in \Sh(p_1,p_2,\ldots,p_k)$ and $L_i\subseteq [p_i]$, $i=2,\ldots, k$. It follows from Lemma \ref{lem:convexlat} and Proposition \ref{prop:convexuvupq} that the statement holds for $k=2$. Now assume for $k-1$.
Write $p=p_1+p_2+\cdots+p_{k-1}$.
Then, by Eq.\,\eqref{eq:shpqp1pk}, there exists $(\eta,\zeta')\in \Sh(p,p_k)\times \Sh(p_1,p_2,\ldots,p_{k-1})$ such that
$\xi=\eta (\zeta'\times 1_{p_k})$. So for any $w=u^{(1)}\times u^{(2)}\times\cdots\times u^{(k)}\in \mathfrak{B}_{p_1}\times  \mathfrak{B}_{p_2,L_2}\times\cdots\times  \mathfrak{B}_{p_k,L_k}$, we have
\begin{align*}
    w \xi^{-1}=\left((u^{(1)}\times\cdots\times u^{(k-1)}) \zeta'^{-1}\times u^{(k)}\right) \eta^{-1}.
\end{align*}
From Proposition \ref{prop:convexuvupq} it follows by induction that $\tau_{\xi}$ is a convex embedding when restricted to
$\mathfrak{B}_{p_1}\times  \mathfrak{B}_{p_2,L_2}\times\cdots\times  \mathfrak{B}_{p_k,L_k}$.
Thus, the restriction of $\tau_{\xi}$ preserves meets and joins by Lemma \ref{lem:convexlat}.
\end{proof}

From Theorem \ref{thm:convexuvupqg}\eqref{lemma-item:abpqccomg} we see that different components
of $\mathfrak{B}_{p_1}\times \mathfrak{B}_{p_2}\times \cdots\times  \mathfrak{B}_{p_k}$ are not adjacent in the Hasse diagram of $\mathfrak{B}_{p_1+p_2 +\ldots+p_k}$ under the weak order.

\begin{exmp}
In the hyperoctahedral group $\mathfrak{B}_3$, as illustrated in Figure $1$, the subgroup $\mathfrak{B}_2\times \mathfrak{B}_1$ is the disjoint union of
\begin{align*}
    \mathfrak{B}_2\times \mathfrak{B}_{1,\emptyset}=[123,\bar1\,\bar23]\quad \text{and}\quad  \mathfrak{B}_2\times \mathfrak{B}_{1,\{1\}}=[12\bar{3},\bar1\,\bar2\,\bar3],
\end{align*}
while  the subgroup $\mathfrak{B}_1\times \mathfrak{B}_2$ is the disjoint union of
\begin{align*}
    \mathfrak{B}_1\times \mathfrak{B}_{2,\emptyset}=[123,\bar1\,32],\quad  \mathfrak{B}_1\times \mathfrak{B}_{2,\{1\}}=[1\bar{2}3,\bar1\,\bar32],\quad
    \mathfrak{B}_1\times \mathfrak{B}_{2,\{2\}}=[13\bar{2},\bar12\bar{3}]
\end{align*}
and
\begin{align*}
    \mathfrak{B}_1\times \mathfrak{B}_{2,\{1,2\}}=[1\bar{3}\,\bar{2},\bar1\,\bar2\,\bar3].
\end{align*}
For signed permutations $132\in \mathfrak{B}_1\times \mathfrak{B}_{2,\emptyset}$ and $13\bar{2}\in \mathfrak{B}_1\times \mathfrak{B}_{2,\{2\}}$,
we have $132<13\bar{2}$ and $\ell(13\bar{2})-\ell(132)=3$. So $3$ is the best lower bound of $\ell(w')-\ell(w)$ when $w'>w$ and they come from different components.

Note that $\Sh(2,1)=\{123,132,231\}$ and $\Sh(1,2)=\{123,213,312\}$. We have
\begin{align*}
    \tau_{132}\left(\mathfrak{B}_2\times \mathfrak{B}_{1,\emptyset}\right)=&\,[132,\bar{1}3\bar{2}],&
    \tau_{231}\left(\mathfrak{B}_2\times \mathfrak{B}_{1,\emptyset}\right)=&\,[312,3\bar{1}\,\bar{2}],\\
    \tau_{132}\left(\mathfrak{B}_2\times \mathfrak{B}_{1,\{1\}}\right)=&\,[1\bar32,\bar1\,\bar{3}\,\bar2],&
    \tau_{231}\left(\mathfrak{B}_2\times \mathfrak{B}_{1,\{1\}}\right)=&\,[\bar{3}12,\bar{3}\,\bar1\,\bar2],\\
    \tau_{213}\left(\mathfrak{B}_1\times \mathfrak{B}_{2,\emptyset}\right)=&\,[213,3\bar{1}2],&
    \tau_{312}\left(\mathfrak{B}_1\times \mathfrak{B}_{2,\emptyset}\right)=&\,[231,32\bar{1}],\\
    \tau_{213}\left(\mathfrak{B}_1\times \mathfrak{B}_{2,\{1\}}\right)=&\,[\bar{2}13,\bar{3}\,\bar{1}2],&
    \tau_{312}\left(\mathfrak{B}_1\times \mathfrak{B}_{2,\{1\}}\right)=&\,[\bar{2}31,\bar{3}2\bar{1}],\\
    \tau_{213}\left(\mathfrak{B}_1\times \mathfrak{B}_{2,\{2\}}\right)=&\,[31\bar{2},2\bar{1}\,\bar{3}],&
    \tau_{312}\left(\mathfrak{B}_1\times \mathfrak{B}_{2,\{2\}}\right)=&\,[3\bar{2}1,2\bar{3}\,\bar{1}],\\
    \tau_{213}\left(\mathfrak{B}_1\times \mathfrak{B}_{2,\{1,2\}}\right)=&\,[\bar{3}1\bar{2},\bar{2}\,\bar{1}\,\bar{3}],&
    \tau_{312}\left(\mathfrak{B}_1\times \mathfrak{B}_{2,\{1,2\}}\right)=&\,[\bar{3}\,\bar{2}1,\bar{2}\,\bar{3}\,\bar{1}].
\end{align*}
\end{exmp}

\section{Monomial basis for the Hopf algebra of signed permutations}\label{sec:mbhsp}

In this section, we apply our results obtained in Sections \ref{sec:lwoohg} and \ref{sec:opmonC} to study the monomial basis $\{\mathbf{M}_{u}\,|\,u\in \mathfrak{B}_n,n\geq0\}$ for the Hopf algebra $\mathfrak{H}Sym$ of signed permutations
with respect to the weak order on hyperoctahedral groups.
In general the structure constants in terms of this basis is not nonnegative.
However, it is shown that for $u\in\mathfrak{B}_{p,[p]}$ the coproduct of $\mathbf{M}_{u}$ is obtained by splitting the signed permutation $u$ at global descents,
and that for $(u,v)\in \mathfrak{B}_{p}\times \mathfrak{B}_{q,[q]}$, the product $\mathbf{M}_{u}\mathbf{M}_{v}$ is a nonnegative linear combination of the monomial basis.
Moreover, under the descent map from $\mathfrak{H}Sym$ to the algebra of type $B$ quasi-symmetric functions,
we show that the image of $\mathbf{M}_{u}$ is either $0$ or a monomial type $B$ quasi-symmetric function.

\subsection{Hopf algebras of permutations and of signed permutations}
Let $\bfk$ be a field of characteristic zero, and let
$$\mathfrak{H}Sym:=\bigoplus_{n=0}^\infty \bfk \mathfrak{B}_n$$
be the graded $\bfk$-vector space whose $n$-th graded component has a linear basis $\mathfrak{B}_n$. By convention $\mathfrak{B}_0$ is the set containing the empty permutation.
Write $\mathbf{F}_{u}$ for the basis element corresponding to the signed permutation $u$ in $ \mathfrak{B}_n$ for $n>0$ and $\iota$ for the basis element of degree $0$.
We call $\{\mathbf{F}_{u}\,|\,u\in \mathfrak{B}_n,n\geq0\}$ the \emph{fundamental basis}.

The space $\mathfrak{H}Sym$ has a self-dual connected graded Hopf algebra structure, whose product and coproduct are induced by the non-parabolic embedding
$\mathfrak{B}_p\times  \mathfrak{B}_{q}\hookrightarrow \mathfrak{B}_{p+q}$.
For $u\in \mathfrak{B}_p$ and $v\in \mathfrak{B}_q$, the product and coproduct are respectively defined by
\begin{align}\label{eq:pdcopdofhsym}
\mathbf{F}_u\mathbf{F}_v=\sum_{\xi\in\Sh(p,q)}\mathbf{F}_{(u\times v)\xi^{-1}}\quad \text{and}\quad
\Delta(\mathbf{F}_u)=\sum_{i=0}^p\mathbf{F}_{\sts(u_1\cdots u_i)}\otimes \mathbf{F}_{\sts(u_{i+1}\cdots u_p)},
\end{align}
where the element $\iota$ is the multiplicative identity.
For instance,
\begin{align*}
\mathbf{F}_{1\bar{2}}\mathbf{F}_{\bar{2}1}=\mathbf{F}_{1\bar{2}\,\bar{4}3}+\mathbf{F}_{1\bar{4}\,\bar{2}3}+\mathbf{F}_{1\bar{4}3\bar{2}}
+\mathbf{F}_{\bar{4}1\bar{2}3}+\mathbf{F}_{\bar{4}13\bar{2}}+\mathbf{F}_{\bar{4}31\bar{2}}
\end{align*}
and
\begin{align*}
\Delta(\mathbf{F}_{1\bar{4}\,\bar{2}3})=\iota\otimes \mathbf{F}_{1\bar{4}\,\bar{2}3}+\mathbf{F}_{1}\otimes \mathbf{F}_{\bar{3}\,\bar{1}2}
+\mathbf{F}_{1\bar{2}}\otimes \mathbf{F}_{\bar{1}2}+\mathbf{F}_{1\bar{3}\,\bar{2}}\otimes \mathbf{F}_{1}+\mathbf{F}_{1\bar{4}\,\bar{2}3}\otimes \iota.
\end{align*}

Let $\mathfrak{S}Sym$ be the vector subspace of $\mathfrak{H}Sym$ generated by  $\{F_{u}\,|\,u\in \mathfrak{S}_n\}$ indexed by permutations, that is
$$
\mathfrak{S}Sym:=\bigoplus_{n=0}^\infty \bfk \mathfrak{S}_n.
$$
Then $\mathfrak{S}Sym$ is indeed a Hopf subalgebra of $\mathfrak{H}Sym$, called the \emph{Malvenuto-Reutenauer Hopf algebra} of permutations.
The product on $\mathfrak{S}Sym$ is defined by Eq.\,\eqref{eq:pdcopdofhsym}, which corresponds to the parabolic embedding
$\mathfrak{S}_p\times  \mathfrak{S}_{q}\hookrightarrow \mathfrak{S}_{p+q}$ of symmetric groups.
Since $\sts$ coincides with $\st$ on permutations, we can write the comultiplication on $\mathfrak{S}Sym$ as
\begin{align*}
\Delta(F_u)=\sum_{i=0}^pF_{\st(u_1\cdots u_i)}\otimes F_{\st(u_{i+1}\cdots u_p)}
\end{align*}
where $u\in \mathfrak{S}_p$.
To illustrate the product and coproduct rules for $\mathfrak{S}Sym$, we have
\begin{align*}
F_{12}F_{21}=F_{1243}+F_{1423}+F_{1432}+F_{4123}+F_{4132}+F_{4312}
\end{align*}
and
\begin{align*}
\Delta(F_{1423})=\iota\otimes F_{1423}+F_{1}\otimes F_{312}+F_{12}\otimes F_{12}+F_{132}\otimes F_{1}+F_{1423}\otimes \iota.
\end{align*}

The Hopf algebra of  permutations was introduced by Malvenuto and Reutenauer \cite{MR95} and extensively studied in \cite{AS05,AT21,DHT02,Pi15,Re05}.
The Hopf algebra $\mathfrak{H}Sym$ was introduced by Aguiar, Bergeron and Nyman \cite{ABN04} in the dual version, and generalized by Novelli and Thibon \cite{NT10} to a family of graded
Hopf algebras with bases labeled by colored permutations. Recently, a non-graded Hopf structure was imposed on the underlying space $\mathfrak{H}Sym$ in \cite{GTY20}, where
it was shown that  $\mathfrak{S}Sym$ is a quotient Hopf algebra of $\mathfrak{H}Sym$.
The nonparabolic embedding $\mathfrak{B}_p\times  \mathfrak{B}_{q}\hookrightarrow \mathfrak{B}_{p+q}$
was used by Huang in \cite{Hu17} to study the representation theory of the hyperoctahedral group, and it is shown that the underlying space $\mathfrak{H}Sym$ has module and comodule structures over the Malvenuto-Reutenauer Hopf algebra $\mathfrak{S}Sym$.

\subsection{Monomial basis}\label{subsec:MonBasis}

Aguiar and Sottile \cite{AS05} defined the monomial basis $\{M_{u}\}$ for $\mathfrak{S}Sym$ with respect to the weak order on $\mathfrak{S}_n$ by
\begin{align*}
    M_{u}:=\sum_{u\leq v}\mu_{\mathfrak{S}_n}(u,v)F_{v},
\end{align*}
where $\mu_{\mathfrak{S}_n}$ is the M\"obius function of the weak order on $\mathfrak{S}_n$.
This basis has a simple coproduct formula.
For a signed permutation $u\in\mathfrak{B}_n$, its \emph{descent set} is
\begin{align*}
    \Des(u):=\{i\in [0,n-1]\,|\,u_i>u_{i+1}\},
\end{align*}
while its \emph{global descent set} is
\begin{align*}
    \GDes(u):=\{i\in [n-1]\,|\,u_j>u_{k}\ \text{for all} j,k\ \text{with}\ 1\leq j\leq i<k\leq n\}.
\end{align*}
Here we use the notation $u_0=0$. Then
\begin{align*}
    \Delta(M_{u})=\sum_{p\in\bar{\GDes(u)}}M_{\st(u_1\cdots u_p)}\otimes M_{\st(u_{p+1}\cdots u_n)},
\end{align*}
where $u\in\mathfrak{S}_n$ and $\bar{\GDes(u)}=\GDes(u)\cup\{0,n\}$.

In an analogous manner, we define a \emph{monomial} basis $\{\mathbf{M}_{u}\}$ for $\mathfrak{H}Sym$ indexed by signed permutations with respect to the weak order on
the hyperoctahedral group.
For each $n\geq0$ and $u\in \mathfrak{B}_n$, let
\begin{align}\label{eq:fundtomono}
    \mathbf{M}_{u}:=\sum_{u\leq  v}\mu_{\mathfrak{B}_n}(u,v)\mathbf{F}_v, \quad \text{or equivalently,}\quad \mathbf{F}_{u}=\sum_{u\leq v}\mathbf{M}_{v}.
\end{align}
where $\mu_{\mathfrak{B}_n}$ is the M\"obius function of the weak order on $\mathfrak{B}_n$.

Although $\mathfrak{S}Sym$ is a Hopf subalgebra of $\mathfrak{H}Sym$ and $F_u=\mathbf{F}_u$ for all permutations $u$, the monomial bases elements $M_u\in\mathfrak{S}Sym$ and $\mathbf{M}_u\in\mathfrak{H}Sym$ indexed by the same permutation $u$ are distinct.
For example,
$$
M_{231}=F_{231}-F_{321},\quad \text{while} \quad\mathbf{M}_{231}=\mathbf{F}_{231}-\mathbf{F}_{321}-\mathbf{F}_{23\bar{1}}+\mathbf{F}_{32\bar{1}}.
$$
So we write $M$ and $F$ in bold type for $\mathfrak{H}Sym$.

\begin{remark}
Recently, Bergeron, D'le\'on, Li, Pang and Vargas \cite{BDLPV21} used an axiomatized version of the methods of Aguiar and Sottile \cite{AS05} on $\mathfrak{S}Sym$
to define a monomial basis on any combinatorial Hopf algebra, and showed that the product is nonnegative and the coproduct is cofree on the monomial basis elements.
It is not difficult to verify that the Hopf algebra $\HSym$ under the weak order on
the hyperoctahedral group
does not satisfy the axioms established in \cite{BDLPV21}.
The following Example \ref{exmp:sccoproduct} shows that the coproduct structure constants of the monomial basis may be negative.

Let $\varphi:\mathfrak{B}_n\rightarrow\mathfrak{S}_n$ be the map forgetting the signs of signed permutations.
Then it follows from the definition that the map $\mathbf{F}_u\mapsto F_{\varphi(u)}$ is a graded Hopf homomorphism from $\HSym$ to $\SSym$.
Note that the map $w\mapsto (\varphi(w), \Nega(w))$ is a bijection between $\mathfrak{B}_n$ and $\mathfrak{S}_n\times 2^{[n]}$, where $2^{[n]}$ is the power set of $[n]$.
Define a partial order $\leq'$ on $\mathfrak{B}_n$ by
$$
u\leq' v\quad\Leftrightarrow \quad \varphi(u)\leq \varphi(v) \quad\text{and}\quad\Nega(u)\subseteq\Nega(v).
$$
It is straightforward to verify that the Hopf algebra $\HSym$  satisfies the axioms given in \cite{BDLPV21} under the order $\leq'$.
\delete{So we obtain a monomial basis with respect to the order $\leq'$ by Eq.\,\eqref{eq:fundtomono} where $\mu_{\mathfrak{B}_n}$ is the M\"obius function of the order $\leq'$ on $\mathfrak{B}_n$.
There exists another order, say $\leq'$, on the hyperoctahedral groups $\mathfrak{B}_n$ such that the monomial basis with respect to this order is assigned by the map $\varphi$ to the monomial basis of $\SSym$ introduced by Aguiar and Sottile \cite{AS05}.}
So we can obtain a monomial basis, with respect to the order $\leq'$, in the sense of \cite{BDLPV21}.
\end{remark}

\begin{exmp}\label{exmp:sccoproduct}
It can be seen from Figure $1$ that
\begin{align*}
\mathbf{M}_{12\overline{3}}=\mathbf{F}_{12\overline{3}}-\mathbf{F}_{21\overline{3}}-\mathbf{F}_{\overline{1}2\overline{3}}
+\mathbf{F}_{\overline{1}\,\overline{2}\,\overline{3}}\quad\text{and}\quad
\mathbf{M}_{1\bar{3}2}=\mathbf{F}_{1\bar{3}2}-\mathbf{F}_{2\overline{3}1}-\mathbf{F}_{\overline{1}\,\overline{3}2}+\mathbf{F}_{\overline{1}\,\overline{3}\,\overline{2}}.
\end{align*}
Direct computations yield that
\begin{align*}
\Delta(\mathbf{M}_{12\bar{3}})=&\,\iota\otimes \mathbf{M}_{12\bar{3}}+\mathbf{M}_{12}\otimes \mathbf{M}_{\bar{1}}-\mathbf{M}_{\bar{1}}\otimes \mathbf{M}_{1\bar{2}}
+\mathbf{M}_{12\bar{3}}\otimes \iota
\intertext{and}
\Delta(\mathbf{M}_{1\bar{3}2})=&\,\iota\otimes \mathbf{M}_{1\bar{3}2}-\mathbf{M}_{\bar{1}}\otimes \mathbf{M}_{\bar{2}1}
-\mathbf{M}_{\bar{1}\,\bar{2}}\otimes \mathbf{M}_{1}+\mathbf{M}_{1\bar{3}2}\otimes \iota.
\end{align*}
\end{exmp}

However, when all entries of $u$ are negative, the coproduct formula of the monomial basis element $\mathbf{M}_u$ is completely analogous to that for $\SSym$.
For instance,
\begin{align*}
    \Delta(\mathbf{M}_{\bar{3}\,\bar{1}\,\bar{2}\,\bar{4}\,\bar{6}\,\bar{5}})=\iota\otimes \mathbf{M}_{\bar{3}\,\bar{1}\,\bar{2}\,\bar{4}\,\bar{6}\,\bar{5}}
    +\mathbf{M}_{\bar{3}\,\bar{1}\,\bar{2}}\otimes \mathbf{M}_{\bar{1}\,\bar{3}\,\bar{2}}+\mathbf{M}_{\bar{3}\,\bar{1}\,\bar{2}\,\bar{4}}\otimes \mathbf{M}_{\bar{2}\,\bar{1}}
    +\mathbf{M}_{\bar{3}\,\bar{1}\,\bar{2}\,\bar{4}\,\bar{6}\,\bar{5}}\otimes \iota.
\end{align*}

\begin{theorem}\label{thm:coformulaofm}
For any $u\in \mathfrak{B}_{n,[n]}$, we have
$$
\Delta(\mathbf{M}_u)=\sum_{p\in \bar{\GDes(u)}}\mathbf{M}_{\sts(u_1\cdots u_p)}\otimes \mathbf{M}_{\sts(u_{p+1}\cdots u_{n})},
$$
where $\bar{\GDes(u)}=\GDes(u)\cup\{0,n\}$.
\end{theorem}
\begin{proof}
Let $\Delta':\HSym\rightarrow \HSym\otimes\HSym$ be the map defined by the right-hand side sum in the desired identity.
Then it suffices to show that $\Delta'$ coincides with the coproduct $\Delta$ defined by Eq.\,\eqref{eq:pdcopdofhsym}.
For any $w\in \mathfrak{B}_{n,[n]}$ and nonnegative integer $p$ with $0\leq p\leq n$, we write $w_{(1)}^{p}=\sts(w_1\cdots w_p)$ and $w_{(2)}^{p}=\sts(w_{p+1}\cdots w_n)$.
Then $p\in \bar{\GDes(w)}$ is equivalent to $w=w_{(1)}^p\times w_{(2)}^p$. Let $u\in \mathfrak{B}_{n,[n]}$. If $v\in \mathfrak{B}_n$ with $u\leq v$, then by Theorem \ref{thm:leqliffinvnegnsp},
$v\in \mathfrak{B}_{n,[n]}$. Thus,
\begin{align*}
    \Delta'(\mathbf{F}_{u})=\sum_{u\leq v}\Delta'(\mathbf{M}_{v})=\sum_{u\leq v}\sum_{p\in \bar{\GDes(v)}}\mathbf{M}_{v_{(1)}^p}\otimes \mathbf{M}_{v_{(2)}^p}
    =\sum_{p=0}^n\sum_{u\leq v\atop v=v_{(1)}^p\times v_{(2)}^p}\mathbf{M}_{v_{(1)}^p}\otimes \mathbf{M}_{v_{(2)}^p}.
\end{align*}
According to Theorem \ref{thm:leqliffinvnegnsp} and Lemma \ref{lem:stsuv},
 $u\leq v_{(1)}^p\times v_{(2)}^p$ if and only if $u_{(1)}^p\leq v_{(1)}^p$ and $u_{(2)}^p\leq v_{(2)}^p$. Consequently,
\begin{align*}
    \Delta'(\mathbf{F}_{u})=&\sum_{p=0}^n\sum_{u\leq v_{1}\times v_{2}}\mathbf{M}_{v_{1}}\otimes \mathbf{M}_{v_{2}}
    =\sum_{p=0}^n\sum_{u_{(1)}^p\leq v_{1}}\mathbf{M}_{v_{1}}\otimes \sum_{u_{(2)}^p\leq v_{2}}\mathbf{M}_{v_{2}}\\
    =&\sum_{p=0}^n\mathbf{F}_{u_{(1)}^p}\otimes \mathbf{F}_{u_{(2)}^p}=\Delta(\mathbf{F}_{u}),
\end{align*}
as desired.
\end{proof}

In general, the product structure constants of $\HSym$ in terms of its monomial basis are not nonnegative.
For example,
\begin{align*}
    \mathbf{M}_{12}\mathbf{M}_{1}={}&\mathbf{M}_{123}+2\mathbf{M}_{132}+\mathbf{M}_{13\bar2}+\mathbf{M}_{231}+\mathbf{M}_{23\bar1}+\mathbf{M}_{312}-\mathbf{M}_{1\bar32}-\mathbf{M}_{\bar{3}12}\\
\intertext{and}
    \mathbf{M}_{\bar{1}2}\mathbf{M}_{1}={}&\mathbf{M}_{\bar{1}23}+2\mathbf{M}_{\bar{1}32}+\mathbf{M}_{\bar{1}3\bar{2}}+\mathbf{M}_{\bar{2}3\bar{1}}
    +\mathbf{M}_{\bar{2}31}+\mathbf{M}_{3\bar{1}2}-\mathbf{M}_{\bar{1}\,\bar{3}2}-\mathbf{M}_{\bar{3}\,\bar{1}2}.
\end{align*}
However, we will deduce from Proposition \ref{prop:convexuvupq} that the product $\mathbf{M}_{u}\mathbf{M}_{v}$ is a nonnegative linear combination of monomial basis for any $(u,v)\in \mathfrak{B}_p\times\mathfrak{B}_{q,[q]}$. For instance,
\begin{align*}
    \mathbf{M}_{1}\mathbf{M}_{\bar{2}\,\bar1}={}&\mathbf{M}_{1\bar{3}\,\bar2}+\mathbf{M}_{\bar31\bar2}+\mathbf{M}_{\bar{3}\,\bar21}.
\end{align*}

Given $u\in \mathfrak{B}_p$, $v\in\mathfrak{B}_q$ and $w\in \mathfrak{B}_{p+q}$,
define $B_{u,v}^w$ to be the set of $\xi\in \Sh(p,q)$ satisfying
\begin{enumerate}
  \item $(u\times v) \xi^{-1}\leq w$;
  \item if $u\leq u'$ such that $(u'\times v)\xi^{-1}\leq w$, then $u=u'$.
\end{enumerate}
Also, define $C_{u,v}^w$ to be the set of $\xi\in \Sh(p,q)$ satisfying
\begin{enumerate}
  \item $(u\times v) \xi^{-1}\leq w$;
  \item if $u\leq u'$ and $v\leq v'$ such that $(u'\times v') \xi^{-1}\leq w$, then $u=u'$ and $v=v'$.
\end{enumerate}
Denote $b_{u,v}^{w}=\#B_{u,v}^w$ and $c_{u,v}^{w}=\#C_{u,v}^w$. Then we have the following theorem.

\begin{theorem}\label{thm:productformulaofm}
Let $p,q$ be positive integers. Then for any $u\in\mathfrak{B}_p$ and $v\in\mathfrak{B}_{q}$, we have
\begin{align}\label{eq:mumvprod1n}
    \mathbf{M}_{u} \mathbf{M}_{v}=\sum_{w\in\mathfrak{B}_{p+q}}\sum_{v\leq v'}\mu(v,v')b_{u,v'}^w\mathbf{M}_{w}.
\end{align}
Moreover, if $\Nega(v)=[q]$, then
\begin{align}\label{eq:mumvprod1post2}
    \mathbf{M}_{u} \mathbf{M}_{v}=\sum_{w\in\mathfrak{B}_{p+q}}c_{u,v}^w\mathbf{M}_{w}.
\end{align}
\end{theorem}
\begin{proof}
We first consider Eq.\,\eqref{eq:mumvprod1n}.
Applying Eq.\,\eqref{eq:fundtomono} to the product $\mathbf{M}_{u} \mathbf{M}_{v}$ gives that
\begin{align*}
     \mathbf{M}_{u} \mathbf{M}_{v}&=\sum_{u\leq u',v\leq v'}\mu{(u,u')}\mu{(v,v')}F_{u'}F_{v'}
     =\sum_{\xi\in \Sh(p,q)}\sum_{u\leq u',v\leq v'}\mu(u,u')\mu(v,v')F_{(u'\times v') \xi^{-1}}.
\end{align*}
Then expressing the fundamental basis in terms of monomial basis yields that
\begin{align}\label{eq:mucdotmv}
     \mathbf{M}_{u} \mathbf{M}_{v}&=\sum_{\xi\in \Sh(p,q)}\sum_{u\leq u', v\leq v'\atop (u'\times v') \xi^{-1}\leq w}\mu(u,u')\mu(v,v')\mathbf{M}_{w}
     =\sum_{w\in\mathfrak{B}_{p+q}}\sum_{u\leq u', v\leq v'}\mu(u,u')\mu(v,v')a_{u',v'}^w\mathbf{M}_{w},
\end{align}
where $a_{u',v'}^w$ is the cardinality of the set
\begin{align*}
    A_{u',v'}^w=\{\xi\in\Sh(p,q)\,|\,(u'\times v') \xi^{-1}\leq w\}.
\end{align*}
Thus, if we can show that
$$
A_{u,v'}^w=\biguplus_{u\leq u'}B_{u',v'}^w,
$$
then, by M\"obius inversion on $\mathfrak{B}_p$, we have
$$
b_{u,v'}^w=\sum_{u\leq u'}\mu(u,u')a_{u',v'}^w,
$$
and hence the desired identity follows.

To see this, let $\xi\in B_{u',v'}^w\cap B_{u'',v'}^w$, where $u\leq u'$ and $u\leq u''$. Then we have
$$
(u'\times v') \xi^{-1}\leq w \quad \text{and}\quad (u''\times v') \xi^{-1}\leq w.
$$
By Proposition \ref{prop:convexuvupq}, $((u'\vee u'')\times v') \xi^{-1}\leq w$, and hence
$u'=u'\vee u''=u''$. So the union is disjoint.

Let $\xi\in B_{u',v'}^w$ for some $u'>u$. Then it follows from Corollary \ref{coro:uuleqvv} and Proposition \ref{prop:convexuvupq} that $(u\times v') \xi^{-1}\leq (u'\times v') \xi^{-1}\leq w$,
so that $\xi\in A_{u,v'}^w$. Conversely, if $\xi\in A_{u,v'}^w$, then put
$$
u''=\bigvee\{u'\in \mathfrak{B}_p\,|\,u\leq u'\ \text{and}\ (u'\times v') \xi^{-1}\leq w \}.
$$
By Proposition \ref{prop:convexuvupq}, the restriction of $\tau_{\xi}$ to the set $\mathfrak{B}_p\times \mathfrak{B}_{q,\Nega(v')}$ is convex, so $\tau_{\xi}$ is join-preserving and we have $(u''\times v') \xi^{-1}\leq w $. By definition, $u''$ is the maximum element of the set of those $u'$
satisfying $(u'\times v') \xi^{-1}\leq w$, and hence $\xi\in B_{u'',v'}^w$. Consequently, $A_{u,v'}^w=\biguplus_{u\leq u'}B_{u',v'}^w$, and Eq.\,\eqref{eq:mumvprod1n} follows.

We now consider Eq.\,\eqref{eq:mumvprod1post2}. Since $\Nega(v)=[q]$, it follows from Theorem \ref{thm:leqliffinvnegnsp} that
if $v\le v'$ in $\mathfrak{B}_q$, then $v'\in\mathfrak{B}_{q,[q]}$.
By Proposition \ref{prop:convexuvupq}, $\tau_{\xi}$ is a convex embedding when restricted to $\mathfrak{B}_p\times  \mathfrak{B}_{q,[q]}$,
a completely analogous argument shows that
$$
A_{u,v}^w=\biguplus_{u\leq u',v\leq v'}C_{u',v'}^w.
$$
Thus, it follows from M\"obius inversion on $\mathfrak{B}_p\times\mathfrak{B}_{q,[q]}$ that
$$
c_{u,v}^w=\sum_{u\leq u', v\leq v'}\mu(u,u')\mu(v,v')a_{u',v'}^w,
$$
which together with Eq.\,\eqref{eq:mucdotmv} yields Eq.\,\eqref{eq:mumvprod1post2}, completing the proof.
\end{proof}

\subsection{The descent map to quasi-symmetric functions of type $B$}
Chow \cite{Ch01} introduced a type $B$ analogous algebra of quasi-symmetric functions, which admits the fundamental basis $\{F_{\alpha}\}$ and the monomial basis $\{M_\alpha\}$ indexed by pseudo-compositions.

Recall that a \emph{pseudo-composition} of $n$ is a sequence $\alpha=(\alpha_1,\alpha_2,\ldots,\alpha_k)$ of nonnegative integers such that $\alpha_1+\alpha_2+\cdots+\alpha_k=n$ with $\alpha_1\geq0$ and $\alpha_i>0$ for $i\geq 2$.
Let $\PC(n)$ denote the set of pseudo compositions of $n$.
To each pseudo composition $\alpha$, we associate the set
$$D(\alpha):=\{\alpha_1,\alpha_1+\alpha_2,\ldots,\alpha_1+\alpha_2+\cdots+\alpha_{k-1}\}.$$
This gives a bijection between the set $\PC(n)$ and the power set $2^{[0,n-1]}$ of $[0,n-1]$, and refinement of pseudo compositions corresponds to inclusion of subsets,  i.e., $\alpha\leq \beta$ if and only if $D(\alpha)\subseteq D(\beta)$.
So we simply identify the poset $\PC(n)$ of  pseudo compositions with the Boolean lattice $2^{[0,n-1]}$ of subsets of $[0,n-1]$.

Let $X=\{x_0,x_1,x_2,\cdots\}$ be a set of commutating indeterminates.
Given a pseudo-composition $\alpha=(\alpha_1,\alpha_2,\ldots,\alpha_k)$ of $n$, the \emph{monomial}
and \emph{fundamental quasi-symmetric functions of type $B$} corresponding to $\alpha$ are defined by
\begin{align*}
    M_{\alpha}:=\sum_{0< i_2<\cdots<  i_n} x_0^{\alpha_1}x_{i_2}^{\alpha_2}\cdots x_{i_k}^{\alpha_k}  \quad \text{and}\quad F_{\alpha}:=\sum_{{0\leq i_1\leq i_2\leq \cdots\leq  i_n}\atop  {j\in D(\alpha)\Rightarrow i_j<i_{j+1}}}
    x_{i_1}x_{i_2}\cdots x_{i_n},
\end{align*}
respectively, where $i_0=0$. Then
\begin{align*}
    F_{\alpha}=\sum_{\alpha\leq \beta}M_{\beta}  \quad \text{and}\quad M_{\alpha}=\sum_{\alpha\leq \beta}(-1)^{\ell(\beta)-\ell(\alpha)}M_{\beta}.
\end{align*}
Let $\BQSym:=\bigoplus_{n\geq 0}\BQSym_n$ where $\BQSym_n$ denotes the space spanned by $\{M_{\alpha}\,|\,\alpha\in\mathcal{PC}(n)\}$, or equivalently by $\{F_{\alpha}\,|\,\alpha\in\mathcal{PC}(n)\}$.
Then $\BQSym$ is a subalgebra of $\bfk[[X]]$, called the \emph{algebra of quasi-symmetric functions of type $B$}.

Chow \cite[Proposition 2.2.6]{Ch01} showed that
the map
\begin{align*}
  \mathcal{D}:\HSym\rightarrow \BQSym,\qquad\mathcal{D}(\mathbf{F}_u)=F_{\Des(u)}
\end{align*}
is a homomorphism of algebras.
We will describe this map in terms of monomial basis. To this end, we first study a relationship between Galois connections and the transfers of monomial bases
on two spaces generated by posets.
A Galois connection is an important tool used in \cite{AS05} to establish the relationship between the monomial bases of the Malvenuto-Reutenauer Hopf algebra of permutations
and the Hopf algebra of quasi-symmetric functions.

Let $P$ and $Q$ be posets, let $f:P\rightarrow Q$ and $g:Q\rightarrow P$ be order-preserving maps.
The pair $(f,g)$ is said to be a \emph{Galois connection} between $P$ and $Q$ if for any $x\in P$
and $y\in Q$,
$$
f(x)\leq_Q y\Leftrightarrow x\leq_P g(y).
$$
When this occurs, Rota \cite[Theorem 1]{Ro64} showed that the M\"obius functions of $P$ and $Q$ are related by
\begin{align}\label{eq:Rotathm}
    \sum_{p\leq x, f(x)=q}\mu_{P}(p,x)=\sum_{y\leq q,\,g(y)=p}\mu_{Q}(y,q)
\end{align}
for all $p\in P$ and $q\in Q$. A generalization of this formula
in the context of Hopf algebras can be found in \cite{AS00}.

When $f$ is surjective we provide a necessary and sufficient condition for a pair of order-preserving maps $(f,g):P\rightleftarrows Q$ to be a Galois connection, which may be useful in its own right.
\begin{lemma}\label{lem:Galoisconcfsurj}
Let $P$ and $Q$ be posets, and let $f:P\rightarrow Q$ and $g:Q\rightarrow P$ be order-preserving maps, where $f$ is surjective.
Then $(f,g):P\leftrightarrows Q$ is a Galois connection if and only if for any $y\in Q$,
$g(y)$ is the largest element of $f^{-1}(y):=\{x\in P\,|\, f(x)=y\}$. In this case, $fg=\Id_Q$.
\end{lemma}
\begin{proof}
If $(f,g)$ is a Galois connection between $P$ and $Q$, then for any $y\in Q$ and any $x\in f^{-1}(y)$,
we have $f(x)\leq y$ so that $x\leq g(y)$.
Thus, $g(y)$ is an upper bound of the set $f^{-1}(y)$.
Since $f$ is order-preserving, we have $y=f(x)\leq f(g(y))$.
On the other hand, it follows from $g(y)\leq g(y)$ that $f(g(y))\leq y$, and hence $ f(g(y))=y$, so $g(y)\in f^{-1}(y)$. Consequently, $g(y)$ is the largest element of  the set $f^{-1}(y)$.

Conversely, assume that $g(y)$ is the largest element of $f^{-1}(y)$ for any $y\in Q$. Then $fg$ is the identity map $\Id_Q$.
Moreover, for any $x\in P$, $g(f(x))$ is the largest element of $f^{-1}(f(x))$, so that $x\leq g(f(x))$. Take any $y\in Q$. If $f(x)\leq y$, then
$x\leq g(f(x))\leq g(y)$ because $g$ is order-preserving. If $x\leq g(y)$, then $f(x)\leq f(g(y))=\Id_Q(y)=y$ because $f$ is order-preserving.
Thus, $(f,g):P\leftrightarrows Q$ is a Galois connection.
\end{proof}

Given a poset $P$, let $\{F_p\,|\,p\in P\}$ be the starting ``fundamental" basis for the space $\bfk P$ and define the ``monomial" basis
 $\{M_p\,|\,p\in P\}$ by
\begin{align*}
    M_{p}:=\sum_{p\leq p'}\mu_{P}(p,p')F_{p'}, \quad \text{or equivalently,}\quad F_p:=\sum_{p\leq p'}M_{p'}.
\end{align*}
Then we have the following result.

\begin{prop}\label{prop:pJGctwo}
Let $P$ and $Q$ be posets with a Galois connection $(f,g):P\rightleftarrows Q$ where $f$ is a surjection.
If $\pi:\bfk P\rightarrow\bfk Q$ is a linear map induced by $\pi(F_p)=F_{f(p)}$ where $p\in P$, then
\begin{align*}
    \pi(M_{p})=\begin{cases}
M_{y},&y\in Q\ \text{and}\ p=g(y),\\
0,&otherwise.
\end{cases}
\end{align*}
\end{prop}
\begin{proof}
It follows from $M_{p}=\sum_{p\leq  x}\mu_{P}(p,x)F_x$ that
\begin{align*}
    \pi(M_p)=\sum_{p\leq x}\mu_{P}(p,x)F_{f(x)}=\sum_{q\in Q}\Big( \sum_{p\leq x,\ f(x)=q}\mu_{P}(p,x)\Big)F_q.
\end{align*}
Since $(f,g):P\rightleftarrows Q$ is a Galois connection,
we see by Eq.\,\eqref{eq:Rotathm} that
\begin{align*}
    \pi(M_p)=\sum_{q\in Q}\Big(\sum_{y\leq q,\,g(y)=p}\mu_{Q}(y,q)\Big)F_q.
\end{align*}
Note that $f$ is a surjection, so by Lemma \ref{lem:Galoisconcfsurj}, we have $fg=\Id_Q$ so that $g$ is injective.
Then there exists a unique element in $Q$, say $y_p$, satisfying $p=g(y_p)$ if $p\in \rm{Im}\,g$; otherwise, the index set $\{y\in Q\,|\,y\leq q,\,g(y)=p\}$ is empty.
Thus,
\begin{align*}
    \sum_{y\leq q,\,g(y)=p}\mu_{Q}(y,q)=\begin{cases}
\mu_{Q}(y,q),&y\in Q\ and\ p=g(y),\\
0,&otherwise.
\end{cases}
\end{align*}
Consequently, $\pi(M_{p})=M_{y}$ if $y$ is the (unique) element in $Q$ such that $p=g(y)$, and  $\pi(M_{p})=0$ otherwise.
\end{proof}

Now we are ready to establish a Galois connection between $\mathfrak{B}_n$ and $2^{[0,n-1]}$,
so that we can use Proposition \ref{prop:pJGctwo} to obtain the relationship between the monomial bases of $\HSym$ and $\BQSym$.

Pseudo compositions can encode descent classes of signed permutations.
Given $I\subseteq [0,n-1]$, let $Y_I$ denote the descent class of $I$ in $\mathfrak{B}_n$, that is, $Y_I:=\{w\in B_n\,|\,\Des(w)=I\}$.
It is well known (see \cite[Theorem 6.2]{BW88}) that $Y_I$ is an interval of $\mathfrak{B}_n$.
There is a concrete description of the largest element of $Y_I$.

\begin{lemma}\label{lem:intdesI}
Let $I$ be a subset of $[0,n-1]$ and $\zeta_I$ the largest element of $Y_I$.
\begin{enumerate}
\item\label{item:maxd0} If $I=\{0<p_1<p_2<\cdots<p_k\}$, then
$$\zeta_I=(\overline{p_1},\ldots, \overline{1},\overline{p_2},\ldots,\overline{p_1+1},\ldots,\overline{n},\ldots,\overline{n-p_k}).$$
\item\label{item:maxd1} If $I=\{p_1<p_2<\cdots<p_k\}$ where $p_1>0$, then
$$\zeta_I=(1,\ldots, p_1,\overline{p_2},\ldots,\overline{p_1+1},\ldots,\overline{n},\ldots,\overline{n-p_k}).$$
\end{enumerate}
\end{lemma}
\begin{proof}
It is straightforward to verify that $\Des(\zeta_I)=I$, so by Theorem \ref{thm:leqliffinvnegnsp} it suffices to show that if $w\in Y_I$, then
$\Nega(w)\subseteq \Nega(\zeta_{I})$, $\Inv(w)\subseteq \Inv(\zeta_{I})$ and $\Nsp(w)\subseteq \Nsp(\zeta_{I})$.
Let $w$ be any signed permutation of $Y_I$. Then $w_{p_1}>w_{p_1+1}, w_{p_2}>w_{p_2+1},\ldots, w_{p_k}>w_{p_k+1}$ and
\begin{align*}
w_1<\cdots< w_{p_1},w_{p_1+1}<\cdots<w_{p_2},\ldots, w_{p_{k}+1}<\cdots<w_n.
\end{align*}
If $I=\{0<p_1<p_2<\cdots<p_k\}$, then $\Nega(\zeta_I)=[n]$, $\Nsp(\zeta_I)=\{(i,j) \,|\,1\leq i<j\leq n\}$ and
\begin{align*}
\Inv(\zeta_I)=\{(i,j)\,|\, 1\leq i\leq p_r<j\leq n\ \text{for some}\ r\in[k]\}.
\end{align*}
If  $I=\{p_1<p_2<\cdots<p_k\}$ where $p_1>0$, then $\Nega(\zeta_I)=[p_1+1,n]$,
\begin{align*}
\Nsp(\zeta_I)=\{(i,j) \,|\,1\leq i\leq p_1<j\leq n\}\cup\{(i,j) \,|\,p_1+1\leq i<j\leq n\}
\intertext{and}
\Inv(\zeta_I)=\{(i,j)\,|\, 1\leq i\leq p_r<j\leq n\ \text{for some}\ r\in[k]\}.
\end{align*}
Moreover, in the former case, we have $w_1<0$, while in the latter case, we have $w_1>0$.
Consequently, in any cases, we must have
$\Nega(w)\subseteq \Nega(\zeta_{I})$, $\Inv(w)\subseteq \Inv(\zeta_{I})$ and $\Nsp(w)\subseteq \Nsp(\zeta_{I})$, and hence $w\leq \zeta_I$.
\end{proof}

Combining Proposition \ref{prop:pJGctwo} and Lemma \ref{lem:intdesI}, the map $\mathcal{D}:\HSym\rightarrow \BQSym$ can be described in terms of the monomial basis.

\begin{theorem}\label{thm:HHSymtoBQSymmonomialb}
For any $w\in\mathfrak{B}_n$, we have
\begin{align*}
    \mathcal{D}(\mathbf{M}_w)=\begin{cases}
    M_{\Des(w)} & \text{if}\ w=\zeta_{\Des(w)},\\
    0& \text{otherwise.}
    \end{cases}
\end{align*}
\end{theorem}
\begin{proof}
By Proposition \ref{prop:pJGctwo}, it suffices to show that the pair of maps $(\Des,g):\mathfrak{B}_n\leftrightarrows 2^{[0,n-1]}$ is a Galois connection,
where $g: 2^{[0,n-1]}\rightarrow \mathfrak{B}_n$ is the map defined by $g(I)=\zeta_I$.
Observe that for any signed permutation $u$, we have $p\in \Des(u)$ is equivalent to $(p,p+1)\in \Inv(u)$.
It then follows from Theorem \ref{thm:leqliffinvnegnsp} that the map $\Des$ is order preserving. According to Lemma \ref{lem:intdesI},
\begin{align*}
g(I)=\zeta_I=\max\{w\in\mathfrak{B}_n\,|\,\Des(w)=I\}=\max\{w\in\mathfrak{B}_n\,|\,\Des(w)\subseteq I\},
\end{align*}
so $g$ is also order preserving. Note that $\Des$ is surjective, so by Lemma \ref{lem:Galoisconcfsurj}, $(\Des,g):\mathfrak{B}_n\leftrightarrows 2^{[0,n-1]}$ is a Galois connection,
and the proof follows from Proposition \ref{prop:pJGctwo} and Lemma \ref{lem:intdesI}.
\end{proof}

\noindent
{\bf Acknowledgements.}
The author would like to thank Li Guo and Zhengpan Wang for helpful conversations.
This work was partially supported by the National Natural Science Foundation of China (Grant Nos. 12071377 and 12071383).

\end{document}